\documentclass{article}
\usepackage{amsfonts,latexsym}
\usepackage{amsmath, times}
\usepackage{amssymb}
\usepackage{amsthm}

\newcommand{\R}{\mathbb R}

\newtheorem{theorem}{Theorem}[section]
\newtheorem{lemma}{Lemma}
\newtheorem{proposition}{Proposition}[section]

\newtheorem{definition}{Definition}[section]

\newlength{\defbaselineskip}
\newcommand{\setlinespacing}[2]%
          {\setlength{\baselineskip}{#1 \defbaselineskip}}

\makeatletter

\@addtoreset{equation}{section} \makeatother 
\begin{document}

\begin{center}
 {\Large   {  Existence solutions for a weighted equation of p-biharmonic type in the unit ball of $\mathbb{R}^{N}$  with critical exponential growth }}
\end{center}
\vspace{0.2cm}

\begin{center}
  Rached Jaidane

 \

\noindent\footnotesize  Department of Mathematics, Faculty of Science of Tunis, University of Tunis El Manar, Tunisia.\\
 Address e-mail: rachedjaidane@gmail.com\\
\end{center}

\vspace{0.5cm}
\noindent {\bf Abstract.}  We study  a weighted $\frac{N}{2}$ biharmonic  equation involving  a positive continuous potential in $\overline{B}$.
 The non-linearity  is assumed to have critical exponential growth in view of logarithmic weighted Adams' type inequalities in the unit ball of $\mathbb{R}^{N}$.   It is proved that there is a nontrivial weak solution to this problem by the mountain Pass Theorem. We avoid the loss of compactness by proving a concentration compactness result and by a suitable asymptotic condition.  \\
\\

\noindent {\footnotesize\emph{Keywords:}Adams' inequality, Moser-Trudinger's inequality, Nonlinearity of exponential growth, Mountain pass method, Compactness level.\\
\noindent {\bf $2010$ Mathematics Subject classification}: $35$J$20$, $35$J$30$, $35$K$57$, $35$J$60$.}

\section{Introduction}

In  this paper, we deal with the  following weighted fourth order problem
\begin{equation}\label{eq:1.1}
  \displaystyle \left\{
      \begin{array}{rclll}
    L:=\Delta(w(x)|\Delta u|^{\frac{N}{2}-2} \Delta u)-\text{div}(|\nabla u|^{\frac{N}{2}-2}\nabla u)+V(x)|u|^{\frac{N}{2}-2}u &=&  \displaystyle f(x,u)& \mbox{in} & B \\
        u=\frac{\partial u}{\partial n}&=&0 &\mbox{on }&  \partial B,
      \end{array}
    \right.
\end{equation}
where the weight $w(x)$ is given by \begin{equation}\label{eq: 1.2}
w(x)=(\log\frac{e}{\vert x\vert})^{\beta(\frac{N}{2}-1)},  \beta\in (0,1),
\end{equation}$B$ is the unitary disk in $\mathbb{R}^{N}$, $N\geq3$, $f(x,t)$ is continuous in
$B\times \mathbb{R}$ and behaves like $\exp\{\alpha t^{\frac{N}{(N-2)(1-\beta)}}\}$ as $ |t|\rightarrow+ \infty$, for some $\alpha >0$ uniformly with respect to $x\in B$.   The potential $V :\overline{B}\rightarrow \mathbb{R}$ is a positive continuous function and bounded away from zero in $B$.
\\

We will look at the historical origins of Trudinger-Moser and Adams inequalities. Trudinger-Moser inequalities are closely linked to the critical case of Sobolev embedding. Trudinger-Moser type inequalities have been established by S. Pohozaev
\cite{Po}, N. Trudinger \cite{NST} and V. Yudovich \cite{Ya}. They essentially showed that for a certain small positive value of $\alpha > 0$, the first-order Sobolev space $$ W_{0}^{1,N}(\Omega)=\mbox{closure}\{u\in C_{0}^{\infty}(\Omega)~~|~~\int_{\Omega}|\nabla u|^{^{N}}dx <\infty\}, $$ endowed with the norm \begin{align}\nonumber
\|u\|_{W^{1,N}_{0}(\Omega)}:=|\nabla u|_{N}=\Big(\int_{\Omega}|\nabla u|^{N}~dx\Big)^{\frac{1}{N}},\end{align} is continuously embedded into the Orlicz space $L_{\phi} (\Omega)$, where $\Omega$ is a smooth bounded domain in $\R^{N}$, $N\geq2$ and $ \phi (t) = e^{\vert t \vert^{\frac{N}{N-1}}}.$
Subsequently, J. Moser improved upon this result. More precisely, he showed that for all $ u\in W_{0}^{1,N}(\Omega),$ $$\exp(\alpha \vert u\vert^{\frac{N}{N-1}})\in L^1(\Omega), \quad \forall\, \alpha > 0\: \mbox{and} $$
\begin{equation}\label{eq:Mt}
\sup_{ \|u\|_{W_{0}^{1,N}(\Omega)}\leq 1}
\int_{\Omega}~e^{\alpha|u|^{\frac{N}{N-1} }}dx < C(N)~~~~\Longleftrightarrow~~~~ \alpha\leq \alpha_{N}:=N{\omega}^{\frac{1}{N-1}}_{N-1},
\end{equation}
where $\omega_{N-1}$ is the area of the unit sphere in $\mathbb{R}^{N}$.  The constant $\alpha_{N}$ is sharp in the sense that for $\alpha > \alpha_{N}$
the supremum in \eqref{eq:Mt} is infinite.\\
Equation (\ref{eq:Mt}) has been used to deal with elliptic problems involving exponential growth nonlinearities. For example, we cite the following problems in dimension $N\geq2$
\begin{equation}\nonumber
-\Delta_{N} u=-\mbox{div}(|\nabla u|^{N-2}\nabla u)= f(x,u)~~\mbox{in}~~\Omega\subset \mathbb{R}^{N},
\end{equation}
which  have been studied considerably by Adimurthi \cite{Adi1,Adi},  Figueiredo et al. \cite{FMR}, Lam and Lu \cite{LL3, LL2, LL1, LL}, Miyagaki and Souto \cite{MS} and Zhang and Chen \cite{ZC}.\\

Significant attention was given to weighted inequalities within weighted Sobolev spaces. This particular type of inequality is recognized in mathematical literature as the weighted Trudinger-Moser inequality \cite{CR1, CR3}. Most studies have predominantly focused on radial functions due to the radial nature of the involved weights. This characteristic amplifies the maximum growth of integrability.

In instances where the weight follows a logarithmic pattern, Calanchi and Ruf expanded the Trudinger-Moser inequality. They established the subsequent results within the weighted Sobolev space, denoted as $W_{0,rad}^{1,N}(B,\rho)=\mbox{closure}\{u \in
C_{0,rad}^{\infty}(B)|\int_{B}|\nabla u|^{N}\rho(x)dx <\infty \}$, where $B$ represents the unit ball in $\mathbb{R}^{N}$.
\begin{theorem}\cite{CR2} \label{th1.1}\begin{itemize}\item[$(i)$] ~~Let $\beta\in[0,1)$ and let $\rho$ given by $ \rho(x)=\big(\log \frac{1}{|x|}\big)^{\beta}$, then
$$
 \int_{B} e^{|u|^{\gamma}} dx <+\infty, ~~\forall~~u\in W_{0,rad}^{1,N}(B,\rho),~~
  \mbox{if and only if}~~\gamma\leq \gamma_{N,\beta}=\frac{N}{(N-1)(1-\beta)}=\frac{N'}{1-\beta}
$$
and
 $$
 \sup_{\substack{u\in W_{0,rad}^{1,N}(B,\rho) \\ \int_{B}|\nabla u|^{N}w(x)dx\leq 1}}
 \int_{B}~e^{\alpha|u|^{\gamma_{N,\beta} }}dx < +\infty~~~~\Leftrightarrow~~~~ \alpha\leq \alpha_{N,\beta}=N[\omega^{\frac{1}{N-1}}_{N-1}(1-\beta)]^{\frac{1}{1-\beta}}
$$
where $\omega_{N-1}$ is the area of the unit sphere $S^{N-1}$ in $\R^{N}$ and $N'$ is the H$\ddot{o}$lder conjugate of $N$.
\item [$(ii)$] Let $\rho$ given by $\rho(x)=\big(\log \frac{e}{|x|}\big)^{N-1}$, then
  \begin{equation*}\label{eq:71.5}
 \int_{B}exp\{e^{|u|^{\frac{N}{N-1}}}\}dx <+\infty, ~~~~\forall~~u\in W_{0,rad}^{1,N}(B,\rho)
 \end{equation*} and
 \begin{equation*}\label{eq:71.6}
\sup_{\substack{u\in W_{0,rad}^{1,N}(B,\rho) \\  \|u\|_{\rho}\leq 1}}
 \int_{B}exp\{\beta e^{\omega_{N-1}^{\frac{1}{N-1}}|u|^{\frac{N}{N-1}}}\}dx < +\infty~~~~\Leftrightarrow~~~~ \beta\leq N,
 \end{equation*}
where $\omega_{N-1}$ is the area of the unit sphere $S^{N-1}$ in $\R^{N}$ and $N'$ is the H$\ddot{o}$lder conjugate of $N$.\end{itemize}
\end{theorem}
The theorem (\ref{th1.1}) has made possible the study of weighted elliptic problems of the second order in dimension $N\geq2$. Thus, Calanchi et al. proved the existence of a non-trivial radial solution for an elliptic problem defined on the $\mathbb{R}^2$ unit ball, with nonlinearities having a double exponential growth at infinity. Subsequently, Deng et al. studied the following problem
\begin{equation}\label{pro}
\displaystyle \left\{
\begin{array}{rclll}
-\textmd{div} (\sigma(x)|\nabla u(x)|^{N-2}\nabla u(x) ) &=& \ f(x,u)& \mbox{in} & B \\
u&=&0 &\mbox{on }&  \partial B,
\end{array}
\right.
 \end{equation}
  where $B$ is the unit ball in $\R^N,\; N\geq2$ and the nonlinearity $f(x, u)$ is continuous in $B\times\R$ and has critical growth in the sense of Theorem (\ref{th1.1}).
The authors have proved that there is a non-trivial solution to this problem, using the mountain pass Theorem. Similar results are proven by Chetouane and Jaidane \cite{CJ}, Dridi \cite{DRI1} and Zhang \cite{Z}. Furthermore, problem \eqref{pro}, involving a potential, has been studied by Baraket and Jaidane \cite{BJ}.
Moreover, we mention that Abid et al. \cite{ABJ} have proved the existence of a positive ground state solution for a weighted second-order elliptic problem of Kirchhoff type, with nonlinearities having a double exponential growth at infinity, using minimax techniques combined with Trudinger-Moser inequality.\\
This notion of critical exponential growth was then extended to higher order Sobolev spaces by Adams' \cite{Ada}. More precisely, Adams' proved the following result,  for $m \in \mathbb{N}$ and $\Omega$ an open bounded set of $\mathbb{R}^N$ such that $m<N,$ there exists a positive constant $C_{m,N}$ such that

\begin{equation}\label{eq:1.31}
\displaystyle\sup_{{u\in W_{0}^{m,\frac{N}{m}}(\Omega), |\nabla^m u|_{\frac{N}{m}} \leq 1}}
\int_{\Omega}~e^{\displaystyle \beta_{0} |u|^{\frac{N}{N-m}}}~dx \leq C_{m,N} \vert \Omega\vert,
\end{equation}
where   $ W_{0}^{m,\frac{N}{m}}(\Omega)$  denotes the $m^{th}$-order Sobolev space, $\nabla^m u$ denotes the $m^{th}$-order gradient of $u$, namely
$$\nabla^m u := \left\{ \begin{array}{ll} \displaystyle\Delta^{\frac{m}{2}} u, \qquad \mbox{if}\; m\; \mbox{is even} \\\\\displaystyle
\nabla\Delta^{\frac{m-1}{2}} u , \qquad \mbox{if}\; m \;\mbox{odd}
\end{array}\right. $$
and
$$\beta_{0} = \beta_{0}(m,N) :=  \frac{N}{\omega_{N-1}} \left\{ \begin{array}{ll}\displaystyle\Big[ \frac{\pi^{\frac{N}{2}} 2^m \Gamma (\frac{m}{2}) }{\Gamma (\frac{N-m}{2}) }\Big]^{\frac{N}{N-m}}, \qquad \mbox{if}\; m\; \mbox{is even} \\\\ \displaystyle
\Big[ \frac{\pi^{\frac{N}{2}} 2^m \Gamma (\frac{m+1}{2}) }{\Gamma (\frac{N-m+1}{2}) }\Big]^{\frac{N}{N-m}}, \qquad \mbox{if}\; m \;\mbox{odd}.
\end{array}\right. $$

  In the particular case where $N = 4$ and
$m = 2,$ the inequality \eqref{eq:1.31} takes the form
\begin{equation}\label{eq:1.311}
\displaystyle\sup_{{u\in W_{0}^{2,2}(\Omega), |\Delta u|_{2} \leq 1}}
\int_{\Omega}~e^{\displaystyle 32 \pi^2 |u|^{2}}~dx \leq C \vert \Omega\vert.
\end{equation}

These findings enabled the exploration of fourth-order problems featuring either subcritical or critical exponential nonlinearity along with continuous potential (refer to the works of Chen et al. \cite{Chen} and Sani \cite{FS}).

Additionally, Ruf and Sani \cite{BRFS} established an inequality similar to Adams', incorporating higher derivatives of even order for unbounded domains in $\mathbb{R}^N$.\\
Recently, an extension of Adams inequalities to Sobolev spaces involving logarithmic weights has been achieved. Wang and Zhu \cite{WZ} have recently established the following result.

 \begin{theorem}\cite{WZ} \label{th1.1} ~~Let $\beta\in(0,1)$ and let $\omega_{\beta}=(\log(\frac{e}{|x|}))^{\beta}$, then
\begin{equation}\label{eq:1.3}
 \displaystyle\sup_{{u\in W_{0,rad}^{2,2}(B,\omega_{\beta}), \|u\| \leq 1}}
 \int_{B}~e^{\displaystyle\alpha|u|^{\frac{2}{1-\beta} }}~dx < \infty~~~~\Leftrightarrow~~~~ \alpha\leq \alpha_{\beta}=4[8\pi^{2}(1-\beta)]^{\frac{1}{1-\beta}},
 \end{equation}
 \end{theorem}
  where $W_{0,rad}^{2,2}(B,\omega_{\beta})$ denotes the weighted Sobolev space of radial functions given by $$W_{0,rad}^{2,2}(B,\omega_{\beta})=\mbox{closure}\Big\{u\in
C_{0,rad}^{\infty}(B)~~|~~\displaystyle\int_{B}\omega_{\beta}(x)|\Delta u|^{2}~dx <\infty\Big\},$$ endowed with the norm  $\displaystyle\|u\|_{W_{0,rad}^{2,2}(B,\omega_{\beta})}= \Big(\int_{B}\omega_{\beta}(x)|\Delta u|^{2}~dx\Big)^{\frac{1}{2}}$, $B=B(0,1)$ is the unit open ball in $\R^{4}$.\\
As an application of Theorem \ref{th1.1}, Dridi and Jaidane \cite{DRI2} considered the following problem
\begin{equation}\label{eq:1.1'}
\displaystyle \left\{
\begin{array}{rclll}
\Delta(\omega_{\beta}(x)\Delta u) - \Delta u +V(x)u&=& f(x,u) \mbox{ in }  B \\
u=\frac{\partial u}{\partial n}&=&0\mbox{ on }  \partial B,
\end{array}
\right.
\end{equation}
where $B=B(0,1)$ is the unit open ball in $\R^{4}$,  $f(x,t)$  is continuous in $B \times \mathbb{R}$ and behaves like $e^{\alpha{t^{\frac{2}{1- \beta}}}}\mbox{ as }t\rightarrow+\infty$, for some $\alpha>0$ and  the potential $V$ is positive and continuous on $\overline{B}$ and bounded away from zero in $B$. The authors proved that there is a nontrivial weak solution to the above problem by using mountain Pass Theorem combined with Trudinger-Moser inequality. Also, Jaidane \cite{RJT} used the same techniques to study a Kirchhoff-type biharmonic problem involving nonlinearities with exponential growth in the sense of the theorem \ref{th1.1}.\\

We denotes by $W^{2,\frac{N}{2}}_{0}(B,w)$ the closure of $C^{\infty}_{0}(B)$ with respect to the norm
\begin{align}\nonumber \|u\|_{W^{2,\frac{N}{2}}}=\displaystyle\big(\int_{B}w(x)|\Delta u|^{\frac{N}{2}}dx +\int_{B} |\nabla u|^{\frac{N}{2}}+\int_{B}|u|^{\frac{N}{2}}dx \big)^{\frac{2}{N}}\cdot\end{align} Also, we consider the subspace of  $W^{2,\frac{N}{2}}_{0}(B,w)$ of radial functions namely  $W^{2,\frac{N}{2}}_{0,rad}(B,w)$.

The choice of the weight  and the space $W_{0,rad}^{2,\frac{N}{2}}(B,w)$ are  motivated by the following exponential inequality (see \cite{ZZ}).
\begin{theorem}\cite{ZZ} \label{th1.1} ~~Let $\beta\in(0,1)$ and let $w$ given by (\ref{eq: 1.2}), then

 \begin{equation}\label{eq:1.3}
 \sup_{\substack{u\in W_{0,rad}^{2,\frac{N}{2}}(B,w) \\  \int_{B}w(x)|\Delta u|^{\frac{N}{2}}dx \leq 1}}
 \int_{B}~e^{\displaystyle\alpha|u|^{\frac{N}{(N-2)(1-\beta)} }}dx < +\infty~~~~\Leftrightarrow~~~~ \alpha\leq \alpha_{\beta}=N[(N-2)N V_{N}]^{\frac{2}{(N-2)(1-\beta)}}(1-\beta)^{\frac{1}{(1-\beta)}},
 \end{equation}
\mbox{where} $V_{N}$ \mbox{ is the volume of the unit ball} $B$ \mbox{in} $\mathbb{R}^{N}$.
 \end{theorem}

Let $\gamma=\gamma(N,\beta):=\displaystyle\frac{N}{(N-2)(1-\beta)}$. In view  of inequality (\ref{eq:1.3}), we say that $f$ has critical growth at $+\infty$ if there exists some $\alpha_{0}>0$,
\begin{equation}\label{eq:1.4}
\lim_{|s|\rightarrow +\infty}\frac{|f(x,s)|}{e^{\alpha s^{\gamma}}}=0,~~~\forall~\alpha~~\mbox{such that}~~ \alpha>\alpha_{0} ~~~~
\mbox{and}~~~~\lim_{|s|\rightarrow +\infty}\frac{|f(x,s)|}{e^{\alpha s^{\gamma}}}=+\infty,~~\forall~ \alpha<\alpha_{0} .\\
\end{equation}

The potential $V$ is continuous on $\overline{B}$ and verifies
\begin{description}
\item[$(V_{1})$] $V(x)\geq V_{0}>0 ~~\mbox{in }~~ B~~ \mbox{for some}~~ V_{0}>0$.

\end{description}
To investigate the solvability of the problem (\ref{eq:1.1}), consider the space $E:= W^{2,\frac{N}{2}}_{0,rad}(B,w)$ with the norm.
\begin{equation*}\label{eq:1.6}
  \| u\|=\left(\int_{B}w(x)|\Delta u|^{\frac{N}{2}}dx+\int_{B}|\nabla u|^{\frac{N}{2}}dx +\int_{B}V(x)|u|^{\frac{N}{2}}dx\right )^{\frac{2}{N}}\cdot
\end{equation*}
 \\

Now, we will present our findings. Throughout this paper, we consistently presume that the nonlinearities $f(x,t)$ exhibit critical growth with $\alpha_{0}> 0$ and fulfill the following conditions:
\begin{description}
\item[$(H_{1})$] The non-linearity $f: \overline{B} \times \mathbb{R}\rightarrow\mathbb{R}$ is continuous, radial in $x$ and $f(x,t)=0$~, $\forall t\leq0$.
  \item[$(H_{2})$]$\mbox{There exist } t_{0} > 0,   M > 0\mbox{ such that } 0 < F(x, t)
=\displaystyle\int_{0}^{t}f(x,s)ds\leq M|f(x,t)|,\newline \forall |t|>t_{0},  \forall
x\in  B$.
    \item[$(H_{3})$] $0<F(x,t)\leq \displaystyle\frac{2}{N}f(x,t)t,\forall
t>0, \forall x\in B.$
\item[$(H_{4})$] $~~~\displaystyle\limsup_{t\rightarrow 0}\frac{N F(x,t)}{2t^{\frac{N}{2}}}< \lambda_{1} ~~~~ \mbox{uniformly in}~~ x,$ where $\lambda_{1} $ is the first eigen value of the operator  with Dirichlet boundary condition  that is \begin{equation}\label{eigen}
\lambda_{1}=\displaystyle\inf_{u\in E, u\neq 0}\frac{\|u\|^{\frac{N}{2}}}{\int_{B}|u|^{\frac{N}{2}}dx}\cdot \end{equation}
This eigen value $\lambda_{1}$ exists and the corresponding  eigen function $\phi_{1}$ is positive and belongs to $L^{\infty}(B)$ \cite{DKN}.
\end{description}

We say that $u$ is a solution to problem (\ref{eq:1.1}), if $u$ is a weak solution
in the following sense.
\begin{definition}\label{def1} A function $u$ is called a solution to $ (\ref{eq:1.1}) $ if $u \in E$ and

\begin{equation}\label{eq:1.6}
\int_{B}\big(w(x)~|\Delta u|^{\frac{N}{2}-2}~\Delta u~~\Delta\varphi+|\nabla u|^{\frac{N}{2}-2} \nabla u.\nabla \varphi +V(x)|u|^{\frac{N}{2}-2}u\varphi\big)~dx=\int_{B}f(x,u)~ \varphi~dx,~~~~~~\mbox{for all }~~\varphi \in E.
\end{equation}
\end{definition}
Finding weak solutions to the problem (\ref{eq:1.1}) is equivalent to identifying non-zero critical points of the following functional inside $E$ :
 \begin{equation}\label{energy}
\mathcal{J}(u)=\frac{2}{N}\int_{B}w(x) |\Delta u|^{\frac{N}{2}} dx+\frac{2}{N}\int_{B}|\nabla u|^{\frac{N}{2}}dx+\frac{2}{N}\int_{B}V(x)|u|^{\frac{N}{2}}dx -\int_{B}F(x,u)dx,
 \end{equation}
 where $F(x,u)=\displaystyle\int_{0}^{u}f(x,t)dt$.\\
\\

We prove the following result.

\begin{theorem}\label{th2}Assume that $V$ is continuous and verifies ($V_{1}$).
Assume that the function $f$  has a critical growth at
$+\infty$ and satisfies
the conditions $(H_{1})$, $(H_{2})$, $(H_{3})$ and $(H_{4})$. If in addition $f$ verifies
the asymptotic condition  $$(H_{5})~~~~~~\displaystyle\lim_{t\rightarrow \infty}\frac{f(x,t)t}{e^{\alpha_{0}t^{\gamma}}}\geq
\gamma_{0}~~~~\mbox{ uniformly in}~~ x, ~~\mbox{with}~~~~\gamma_{0}> \displaystyle \frac{(\frac{\alpha_{\beta}}{\alpha_{0}})^{\frac{N}{2\gamma}}}{V_{N}e^{N(1-\log (2e))}},$$ then the
problem  (\ref{eq:1.1}) has a nontrivial solution.
\end{theorem}

In general, the exploration of fourth-order partial differential equations remains an intriguing field. The interest in these equations sparked due to their applications in various areas such as micro-electro-mechanical systems, phase field models for multi-phase systems, thin film theory, surface diffusion on solids, interface dynamics, and flow in Hele-Shaw cells (refer to \cite{D, FW, M}). However, numerous practical applications stem from elliptic problems, such as the investigation of traveling waves in suspension bridges and radar imaging (see, for instance, \cite{AEG, LL}).\\
This paper is structured as follows:

Section 2 covers essential preliminary information about functional spaces.
Section 3 presents several useful lemmas crucial for compactness analysis.
In Section 4, we demonstrate that the energy $\mathcal{J}$ satisfies two specific geometric properties and meets the compactness condition within a specified level.
Finally, in Section 5, we conclude with the proof of the main result.
Throughout this paper, the constant $C$ might vary between different lines, and at times, we index the constants to illustrate their variations.

\section{Weighted Lebesgue and Sobolev Spaces setting}
Let $\Omega$ be a bounded domain in $\mathbb{R}^{N}$, where $N\geq2$, and let $w$ be a nonnegative function belonging to the space $L^{1}(\Omega)$. To address operators involving weights, we introduce functional spaces such as $L^{p}(\Omega,w)$, $W^{m,p}(\Omega,w)$, and $W_{0}^{m,p}(\Omega,w)$ along with their pertinent properties, which will be utilized later on.

Let $S(\Omega)$ represent the set  of all measurable real-valued functions defined on $\Omega$, where two measurable functions are considered identical if they are almost everywhere equal.

Following the definitions provided by Drabek, Kufner, and others in \cite{DKN}, the weighted Lebesgue space $L^{p}(\Omega,w)$ is defined as the set of measurable functions $u:\Omega\rightarrow \mathbb{R}$ satisfying $\int_{\Omega}w(x)|u|^{p}~dx<\infty$, for any real number $1\leq p<\infty$.

This space is equipped with a norm given by
$$\|u\|_{p,w}=\Big(\int_{\Omega}w(x)|u|^{p}~dx\Big)^{\frac{1}{p}}.$$

For $m\geq 2$, let $w$ be a given family of weight functions $w_{\tau}, ~~|\tau|\leq m,$ $w=\{w_{\tau}(x)~~x\in\Omega,~~|\tau|\leq m\}.$\\
In \cite{DKN}, the  corresponding weighted Sobolev space was  defined as
$$ W^{m,p}(\Omega,w)=\{ u \in L^{p}(\Omega), D^{\tau} u \in L^{p} (\Omega)~~  \forall ~~1\leq|\tau|\leq m-1 , D^{\tau} u   \in L^{p}(\Omega,w) ~~  \forall ~~|\tau|=m  \}$$
endowed with the following norm:

\begin{equation*}\label{eq:2.2}
\|u\|_{W^{m,p}(\Omega,w)}=\bigg(\sum_{ |\tau|\leq m-1}\int_{\Omega}|D^{\tau}u|^{p}dx+\displaystyle \sum_{ |\tau|= m}\int_{\Omega}w(x) |D^{\tau}u|^{p}dx\bigg)^{\frac{1}{p}},
\end{equation*}
where $w_{\tau}=1~~\mbox{for all}~~|\tau|< k,$ $w_{\tau}=w~~\mbox{for all}~~|\tau|=k$.\\
If also we suppose that $w(x)\in L^{1}_{loc}(\Omega)$, then $C^{\infty}_{0}(\Omega)$ is a subset of $W^{m,p}(\Omega,w)$ and we can introduce the space $W^{m,p}_{0}(\Omega,w)$
as the closure of $C^{\infty}_{0}(\Omega)$ in $W^{m,p}(\Omega,w).$  Moroever, the injection $$W^{m,p}(\Omega,w)\hookrightarrow W^{m-1,p}(\Omega)~~\mbox{is compact}.$$
Also, $(L^{p}(\Omega,w),\|\cdot\|_{p,w})$ and $(W^{m,p}(\Omega,w),\|\cdot\|_{W^{m,p}(\Omega,w)})$ are separable, reflexive Banach spaces provided that $w(x)^{\frac{-1}{p-1}} \in L^{1}_{loc}(\Omega)$.
Then the space $E$ is a Banach and reflexive space provided $(V_{1})$ is satisfied. Furthermore, the space $E$ is endowed with the norm
\begin{align}\nonumber\| u\|=\left(\int_{B}w(x)|\Delta u|^{\frac{N}{2}}dx +\int_{B} |\nabla u|^{\frac{N}{2}}+\int_{B}V(x)|u|^{\frac{N}{2}}dx\right )^{\frac{N}{2}}\end{align}
which is equivalent to the following norm (see lemma \ref{lem1}) \begin{align}\nonumber\|u\|_{W_{0,rad}^{2,\frac{N}{2}}(B,w)}=\displaystyle\big(\int_{B}w(x)|\Delta u|^{\frac{N}{2}}dx +\int_{B} |\nabla u|^{\frac{N}{2}}+\int_{B}|u|^{\frac{N}{2}}dx \big)^{\frac{2}{N}}\cdot\end{align}

\section{ Preliminaries for the compactness analysis }

This section will present several technical lemmas that we'll use later on. We'll start with the radial lemma.

\begin{lemma}\label{lem1} Assume that $V$ is continuous and verifies ($V_{1}$).

\item[(i)] Let $u$ be a radially symmetric
 function in $C_{0,rad}^{\infty}(B)$. Then, we have\begin{itemize}\item[$(i)$]\cite{ZZ}
 $$|u(x)|\leq \bigg(\displaystyle\frac{N}{\alpha_{\beta}}\big(|\log(\frac{e}{|x|}|-1\big)\bigg)^{\frac{1}{\gamma}}\displaystyle \big(\int_B w(x)|\Delta u|^{\frac{N}{2}}dx \big)^{\frac{2}{N}}\leq \displaystyle \bigg(\displaystyle\frac{N}{\alpha_{\beta}}\big(|\log(\frac{e}{|x|}|-1\big)\bigg)^{\frac{1}{\gamma}}\|u\|\cdot$$
\item[$(ii)$]$\displaystyle\int_{B}e^{|u|^{\gamma}}dx<+\infty,~~\forall u\in .$

\item[(iii)]  The following embedding
is continuous $$E\hookrightarrow L^{q}(B)~~\mbox{for all}~~q\geq \frac{N}{2}.$$
\item[(vi)]$E$ is compactly
embedded in $L^{q}(B)$  for all $q \geq1$.
\end{itemize}
\end{lemma}
\textit{Proof }
$(i)$ see \cite{WZ}\\
$(ii)$ Follows from $(i)$ and by density.\\
  $(iii)$ Since $w(x)\geq1$ and $V_{0}>0$, then the following embedding are continuous   $$ E = W^{2,\frac{N}{2}}_{0,rad}(B,w)\hookrightarrow W^{2,\frac{N}{2}}_{0,rad}(B)\hookrightarrow L^{q}(B)~~\forall q\geq \frac{N}{2}.$$  We also have, by Rellich-Kondrachov, the following compact injection$$ W^{2,\frac{N}{2}}_{0,rad}(B)\hookrightarrow\hookrightarrow L^{q}(B)~~\forall q\geq 1.$$ So, $ E$ is compactly embedded in $L^{q}(B)~~\forall q\geq 1$.

 This concludes the lemma.\hfill $\Box$\\
Second, we give the following useful lemma.
 \begin{lemma}\cite{FMR}\label{lem2}
 Let $\Omega\subset \mathbb{R^{N}}$ be a bounded domain and $f:\overline{\Omega}\times\mathbb{R}$
  a continuous function. Let $(u_{n})_{n}$ be a sequence in $L^{1}(\Omega)$
converging to $u$ in $L^{1}(\Omega)$. Assume that $\displaystyle f(x,u_{n})$ and
$\displaystyle f(x,u)$ are also in $ L^{1}(\Omega)$. If
$$\displaystyle\int_{\Omega}|f(x,u_{n})u_{n}|dx \leq C,$$ \\where $C$
is a positive constant, then $$f(x,u_{n})\rightarrow
f(x,u)~~\mbox{in}~~L^{1}(\Omega).$$
\end{lemma}
In the sequel, we prove a concentration compactness result of Lions type.
\begin{lemma}\label{Lionstype} Let
$(u_{k})_{k}$  be a sequence in $E$. Suppose that, \newline $\|u_{k}\|=1$, $u_{k}\rightharpoonup u$ weakly in $E$, $u_{k}(x)\rightarrow u(x) ~~a.e~x\in B$, $\nabla u_{k}(x)\rightarrow\nabla u(x) ~~a.e~x\in B$,  $\Delta u_{k}(x)\rightarrow\Delta u(x) ~~a.e~x\in B$ and $u\not\equiv 0$. Then
$$\displaystyle\sup_{k}\int_{B}e^{p~\alpha_{\beta}
|u_{k}|^{\gamma}}dx< +\infty,~~\mbox{where}~~ \alpha_{\beta}=N[(N-2)N V_{N}]^{\frac{2}{(N-2)(1-\beta)}}(1-\beta)^{\frac{1}{(1-\beta)}},$$
for all $1<p<U(u)$ where $U(u)$ is given by:
 $$U(u):=\displaystyle \left\{
      \begin{array}{rcll}
&\displaystyle\frac{1}{(1-\|u\|^{\frac{N}{2}})^{\frac{2\gamma}{N}}}& \mbox{ if }\|u\| <1\\
       &+\infty& \mbox{ if } \|u\|=1\\
 \end{array}
    \right.$$
\end{lemma}
\textit{Proof}
For $a,~b \in \R,~~q>1$. If $q'$ its conjugate i.e. $\frac{1}{q}+\frac{1}{q'}=1$
we have, by Young inequality, that
$$e^{a+b}\leq \frac{1}{q}e^{qa}+ \frac{1}{q'}e^{q'b}.$$

Also, we have
\begin{equation}\label {eq:3.1}
(1+a)^{q}\leq (1+\varepsilon) a^{q}+(1-\frac{1}{(1+\varepsilon)^{\frac{1}{q-1}}})^{1-q},~~\forall a\geq0,~~\forall\varepsilon>0~~\forall q>1.
\end{equation}
So, we get
$$
      \begin{array}{rcll}
|u_{k}|^{\gamma}&=& |u_{k}-u+u|^{\gamma}\\
&\leq& (|u_{k}-u|+|u|)^{\gamma}\\
&\leq& (1+\varepsilon)|u_{k}-u|^{\gamma}+\big(1-\frac{1}{(1+\varepsilon)^{\frac{1}{\gamma-1}}}\big)^{1-\gamma}|u|^{\gamma}\\
 \end{array}
   $$
which implies that
 \begin{align*}
\int_{B}e^{p~\alpha_{\beta}
|u_{k}|^{\gamma}}dx &\leq
\frac{1}{q}\int_{B} e^{pq~\alpha_{\beta}
(1+\varepsilon)|u_{k}-u|^{\gamma}}dx\\
&+\displaystyle\frac{1}{q'}\int_{B} e^{pq'~\alpha_{\beta}
(1-\frac{1}{(1+\varepsilon)^{\frac{1}{\gamma-1}}})^{1-\gamma}|u|^{\gamma}}dx,
\end{align*}
for any $p>1$.
From Lemma \ref{lem1} $(ii)$, the last integral is finite.\\ To finish the
proof, we need to prove that for all $p$ such that $1<p<U(u)$,
\begin{equation}\label {eq:3.2}
\sup_{k}\int_{B} e^{pq~\alpha_{\beta}
(1+\varepsilon)|u_{k}-u|^{\gamma}}dx<+\infty,
\end{equation}
 for some $\varepsilon>0$ and $q>1$.\\
In what follows, we assume that
$\|u\|<1$ and in the case that $\|u\|=1$, the proof  is similar.\\
When $$\|u\|<1$$
and
$$p<\displaystyle\frac{1}{(1-\|u\|^{\frac{N}{2}})^{\frac{2\gamma}{N}}},$$
there exists $\nu>0$ such that
$$p(1-\|u\|^{\frac{N}{2}})^{\frac{2\gamma}{N}}(1+\nu)<1.$$
On the other hand,
 by Brezis-Lieb's Lemma \cite{Br} we have
 \begin{equation*}\label {eq:2.4}
 \|u_{k}-u\|^{\frac{N}{2}}=\|u_{k}\|^{\frac{N}{2}}-\|u\|^{\frac{N}{2}}+o(1)~~\mbox{where}~~
  o(1)\rightarrow 0~~ \mbox{as}~~k\rightarrow +\infty.
\end{equation*}
   Then,
   \begin{align}\nonumber\|u_{k}-u\|^{\frac{N}{2}}=1-\|u\|^{\frac{N}{2}}+o(1),\end{align}
and so

$$\displaystyle\lim_{k\rightarrow+\infty}\|u_{k}-u\|^{\gamma}=(1-\|u\|^{\frac{N}{2}})^{\frac{2\gamma}{N}}.$$

 Therefore, for every
$\varepsilon>0$, there exists $k_{\varepsilon}\geq 1$ such that
$$\|u_{k}-u\| ^{\gamma}\leq (1+\varepsilon)(1-\|u\|^{\frac{N}{2}})^{\frac{2\gamma}{N}},~~\forall ~~k\geq k_{\varepsilon}.$$
If we take $q=1+\varepsilon$ with
$\varepsilon=\sqrt[3]{1+\nu}-1$, then  $\forall k\geq k_{\varepsilon}$, we have
 $$pq(1+\varepsilon)\|u_{k}-u\|^{\gamma}\leq 1.$$
Consequently,
$$\begin{array}{rlll}
\displaystyle\int_{B} e^{pq~\alpha_{\beta}
(1+\varepsilon)|u_{k}-u|^{\gamma}}dx&\leq&
\displaystyle\int_{B} e^{
(1+\varepsilon)pq~\alpha_{\beta}(\frac{|u_{k}-u|}{\|u_{k}-u\|})^{\gamma}\|u_{k}-u\|^{\gamma}}dx\\
&\leq&\displaystyle\int_{B} e^{~\alpha_{\beta}(\frac{|u_{k}-u|}{\|u_{k}-u\|})^{\gamma}}dx\\
&\leq &\displaystyle\sup _{\|u\|\leq 1}\displaystyle\int_{B}
e^{~\alpha_{\beta}|u|^{\gamma}}dx <+\infty.
\end{array}
$$
Now, (\ref{eq:3.2})  follows from  (\ref{eq:1.3}).
 This complete the proof and lemma \ref{Lionstype} is proved.
\section{The variational formulation}
Since the reaction term $f$ is of critical exponential growth, there exist positive constants $a$ and  $
C$ such that
\begin{equation}\label {eq:4.1}
|f(x,t)|\leq C e^{a ~t^{\gamma}}, ~~~~~~\forall |t| >t_{1}\cdot
\end{equation}
 So, by using $(H_{1})$, the functional $\mathcal{J}$ given by (\ref{energy}) is $C^{1}$.
\subsection{The mountain pass geometry of the energy}
In the sequel, we prove that the functional $\mathcal{J}$ has a mountain pass
geometry.

\begin{proposition}\label{prg}
Assume that the hypothesis $(H_{1}),(H_{2}),(H_{3}),(H_{4})$, and $(V_{1})$  hold. Then,
 \begin{enumerate}
 \item[(i)] there exist  $\rho,~\beta_{0}>0$  such that $\mathcal{J}(u)\geq \beta_{0}$ for all
 $u\in E$ with $\|u\|=\rho$.
\item[(ii)] Let $\varphi_{1}$ be a normalized eigen function associated to $\lambda_{1}$ in $E$. Then, $\mathcal{J}(t\varphi_{1})\rightarrow -\infty,~~\mbox{as}~~t\rightarrow+\infty$.
 \end{enumerate}
\end{proposition}
 \textit{Proof}~~It follows from the hypothesis $(H_{4})$ that  there exists $t_{2}>0$ and there exists  $\varepsilon_{0} \in (0,1)$  such that
\begin{equation}\label {eq:4.2}
F(x,t)\leq \frac{2}{N}\lambda_{1}(1-\varepsilon_{0})|t|^{\frac{N}{2}},~~~~~~\mbox{for}~~|t|< t_{2}.
\end{equation}
Indeed,  $$\displaystyle\limsup_{t\rightarrow 0}\frac{N F(x,t)}{2t^{\frac{N}{2}}}< \lambda_{1}$$
or $$\displaystyle\inf_{\tau>0}\sup\{\frac{NF(x,t)}{2t^{\frac{N}{2}}},~~0<|t|<\tau\}< \lambda_{1}$$
Since this inequality is strict, then there exists $\varepsilon_{0}>0$ such that
$$\displaystyle\inf_{\tau>0}\sup\{\frac{NF(x,t)}{2t^{\frac{N}{2}}},~~0<|t|<\tau\}< \lambda_{1}-\varepsilon_{0},$$
hence, there exists $t_{2}>0$ such that
$$\sup\{\frac{NF(x,t)}{2t^{\frac{N}{2}}},~~0<|t|<t_{2}\}< \lambda_{1}-\varepsilon_{0}.$$
Hence$$\forall |t|<t_{2}~~~F(x,t)\leq \frac{2}{N}\lambda_{1}(1-\varepsilon_{0})t^{\frac{N}{2}}.$$

From $(H_{3})$ and (\ref{eq:4.1}) and for all $q>2$, there exist a  constant $C>0$ such that
\begin{equation}\label{eq:4.3}
F(x,t)\leq  C |t|^{q} e^{a~t^{\gamma}},~~\forall~|t| > t_{1}.
\end{equation}
So
\begin{equation}\label {eq:4.4}
F(x,t)\leq \frac{2}{N}\lambda_{1}(1-\varepsilon_{0})|t|^{\frac{N}{2}}+C |t|^{q}e^{a~t^{\gamma}},~~~~~~\mbox{for all}~~t\in \R.
\end{equation}
Since
$$\mathcal{J}(u)=\frac{2}{N}\|u\|^{\frac{N}{2}}-\int_{B}F(x,u)dx,$$
we get
$$\mathcal{J}(u)\geq \frac{2}{N}\|u\|^{\frac{N}{2}}- \frac{2}{N}\lambda_{1}(1-\varepsilon_{0})\|u\|^{\frac{2}{N}}_{\frac{N}{2}}-C\int_{B} |u|^{q}e^{a~u^{\gamma}}~dx.$$
But $\lambda_{1}\|u\|_{\frac{N}{2}}^{2}\leq \|u\|^{\frac{N}{2}}$ and from the  H\"{o}lder inequality, we obtain
\begin{equation}\label {eq:4.5}
\mathcal{J}(u)\geq \frac{ 2\varepsilon_{0}}{N}\|u\|^{\frac{N}{2}}-C  (\int_{B}e^{2a~|u|^{\gamma}}dx\Big)^{\frac{1}{2}}\|u\|^{q}_{2q}\cdot
\end{equation}
From the Theorem \ref{th1.1}, if we choose $u\in E$ such that
\begin{equation}\label {eq:4.6}
2a \|u\|^{\gamma}\leq \alpha_{\beta},
\end{equation}
we get
$$\int_{B}e^{2a|u|^{\gamma}}dx=\int_{B}e^{2a\|u\|^{\gamma}(\frac{|u|}{\|u\|})^{\gamma}}dx<+\infty.$$
On the other hand $\|u\|_{2q}\leq C \|u\|$ (Lemma \ref{lem1}),  so
$$\mathcal{J}(u)\geq \frac{2\varepsilon_{0}}{N}\|u\|^{\frac{N}{2}}-C\|u\|^{q},$$
for all $u\in E$ satisfying (\ref{eq:4.6}).
Since $2<q$, we can choose $\rho=\|u\|>0$ as the maximum point of the function $g(\sigma)=\frac{2\varepsilon_{0}}{N} \sigma^{\frac{N}{2}}-C\sigma^{q}$ on the interval $[0,(\frac{\alpha_{\beta}}{2a})^{\frac{1}{\gamma}}]$ and $\beta_{0}=g(\rho)$ , $\mathcal{J}(u) \geq\beta_{0}>0$.\\

(ii)  Let $\phi_{1}\in E\cap L^{\infty}(B)$ be the normalized eigen function associated to the eigen value defined by (\ref{eigen}) ie such that $\|\phi_{1}\|=1$.
We define the function
$$\displaystyle\varphi
(t)=\mathcal{J}(t\phi_{1})=\frac{2t^{\frac{N}{2}}}{N}\|\phi_{1}\|^{\frac{N}{2}}-\int_{B}F(x,t\phi_{1})
dx.$$ Then using $(H_{1})$,  $(H_{2})$ and integrating, we get the
existence of a constant $C>0$ such that $$F(x,t)\geq C
e^{\frac{1}{M}t},~~\forall ~~|t|\geq t_{0}.$$ In particular, for
$p>\frac{N}{2}$, there exists $C$ such that $$F(x,t)\geq
C|t|^{p}-C,~~\forall t\in \mathbb{R},~~x\in B.$$ Hence, $$\varphi
(t)=\mathcal{J}(t\phi_{1})\leq\frac{t^{\frac{N}{2}}}{\frac{N}{2}}\|\phi_{1}\|^{\frac{N}{2}}-|t|^{p}\|\phi_{1}\|_{p}-c_{5}\rightarrow
-\infty,~~\mbox{as}~~t\rightarrow+\infty,$$ and it's easy to
conclude.\\\\
\subsection{The compactness level of the energy}

The primary challenge within the variational approach to the critical growth problem arises from the absence of compactness. Specifically, the global Palais-Smale condition doesn't hold as required. However, a partial Palais-Smale condition is maintained under a specified level. In the subsequent proposition, we pinpoint the initial level of non-compactness within the energy.

\begin{proposition}\label{propPS} Let $\mathcal{J}$ be the energy associated to the problem (\ref{eq:1.1}) defined by (\ref{energy}), and suppose that the conditions $(V_{1})$, $(H_{1})$, $(H_{2})$, $(H_{3})$ and $(H_{4})$ are satisfied. If the function $f(x,t)$ satisfies the condition (\ref{eq:1.4}) for some
$\alpha_{0} >0$, then the functional $\mathcal{J}$ satisfies the Palais-Smale condition $(PS)_{c}$ for any
$$c<\displaystyle\frac{2}{N}(\frac{\alpha_{\beta}}{\alpha_{0}})^{\frac{N}{2\gamma}},$$
where $\alpha_{\beta}=N[(N-2)N V_{N}]^{\frac{2}{(N-2)(1-\beta)}}(1-\beta)^{\frac{1}{(1-\beta)}}.$

\end{proposition}
 {\it Proof}~~Consider a $(PS)_{c}$ sequence ($u_{n}$) in $E$, for some $c\in \R$, that is
\begin{equation}\label{eq:4.7}
\mathcal{J}(u_{n})=\frac{2}{N}\|u_{n}\|^{\frac{N}{2}}-\int_{B}F(x,u_{n})dx \rightarrow c ,~~n\rightarrow +\infty
\end{equation}
and
 \begin{align}\label{eq:4.8}\nonumber
\mathcal{ J}'(u_{n})\varphi &=\Big|\int_{B}\big(w(x)~|\Delta u_{n}|^{\frac{N}{2}-2}~\Delta u_{n}~~\Delta\varphi+|\nabla u_{n}|^{\frac{N}{2}-2} \nabla u_{n}.\nabla \varphi +V(x)|u_{n}|^{\frac{N}{2}-2}u_{n}\varphi\big)~dx-\int_{B}f(x,u_{n})~ \varphi~dx
\Big|\\&\leq \varepsilon_{n}\|\varphi\|,
\end{align}
for all $\varphi \in E$, where $\varepsilon_{n}\rightarrow0$, as $n\rightarrow +\infty$.\\

 It follows from  $(H_{2})$ , that for all
$\varepsilon>0$ there exists $t_{\varepsilon}>0$ such that
\begin{equation}\label{eq:4.9}
 F(x,t)\leq \varepsilon t f(x,t),~~~~\mbox{for all}~~ |t|>t_{\varepsilon}~~\mbox{and uniformly in}~~ x\in B,
 \end{equation}
and so, by (\ref{eq:4.7}), for large enough $n$, there exists a constant $C>0$
$$\displaystyle\frac{2}{N}\|u_{n}\|^{\frac{N}{2}}\leq C+\displaystyle\int_{B}F(x,u_{n})dx,$$
hence
$$
\displaystyle\frac{2}{N}\|u_{n}\|^{\frac{N}{2}}\leq  C +\displaystyle\int_{{|u_{n}|\leq t_{\varepsilon}}}F(x,u_{n})dx+ \varepsilon \int_{B}f(x,u_{n})u_{n}dx $$
and so, from (\ref{eq:4.8}), we get
$$\displaystyle\frac{2}{N}\|u_{n}\|^{\frac{N}{2}}\leq C_{1}+\varepsilon \varepsilon_{n}\|u_{n}\|+ \varepsilon \|u_{n}\|^{\frac{N}{2}},$$
for some constant $ C_{1}>0$.
Since
\begin{equation}\label{eq:4.10}
\displaystyle(\frac{2}{N}-\varepsilon)\|u_{n}\|^{\frac{N}{2}}\leq C_{1}+\varepsilon \varepsilon_{n}\|u_{n}\|,
 \end{equation}
 we deduce that the sequence $(u_{n})$ is bounded in $E$. As consequence, there exists $u\in E$ such that, up to subsequence,
 $u_{n}\rightharpoonup u $ weakly in $E$, $u_{n}\rightarrow u$ strongly in $L^{q}(B)$, for all $ q \geq1$ and $u_{n}(x)\rightarrow u(x)$ a.e. in $B$ . Also, we can follow  \cite{CJ} to prove that  $\nabla u_{n}(x)\rightarrow\nabla u(x) ~\mbox{a.e}~x\in B$ and $\Delta v_{n}(x)\rightarrow\Delta v(x) ~\mbox{a.e}~x\in B $. \\
Furthermore, we have, from (\ref{eq:4.7}) and (\ref{eq:4.8}), that
\begin{equation}\label{eq:4.11}
0<\int_{B}| f(x,u_{n})u_{n}|~dx\leq C,
 \end{equation}
and
 \begin{equation}\label{eq:4.12}
0<\int_{B} F(x,u_{n})~dx\leq C.
 \end{equation}
Since by Lemma 2.1 in \cite {FMR}, we have
\begin{equation}\label{eq:4.13}
f(x,u_{n})\rightarrow f(x,u) ~~\mbox{in}~~L^{1}(B) ~~as~~ n\rightarrow +\infty,
 \end{equation}
then, it follows from $(H_{2})$ and the generalized Lebesgue dominated convergence theorem that
\begin{equation}\label{eq:4.14}
F(x,u_{n})\rightarrow F(x,u) ~~\mbox{in}~~L^{1}(B) ~~as~~ n\rightarrow +\infty.
 \end{equation}
So, using (\ref{eq:4.14}), we get
\begin{equation}\label{eq:4.15}
\displaystyle
\lim_{n\rightarrow+\infty}\|u_{n}\|^{\frac{N}{2}}=\frac{N}{2}(c+\int_{B}F(x,u)dx).
 \end{equation}
Then from (\ref{eq:4.4}), we have
\begin{equation}\label{eq:4.16}
\displaystyle\lim_{n\rightarrow+\infty}\int_{B}f(x,u_{n})u_{n}dx=\frac{N}{2}(c+\int_{B}F(x,u)dx).
 \end{equation}
It follows from $(H_{3})$ and (\ref{eq:4.8}), that
\begin{equation}\label{eq:4.17}
\displaystyle
\lim_{n\rightarrow+\infty}\frac{N}{2} \int_{B}F(x,u_{n})dx\leq \displaystyle
\lim_{n\rightarrow+\infty}\int_{B}f(x,u_{n})u_{n}dx= \frac{N}{2}(c+\int_{B}F(x,u)dx).
\end{equation}
As a direct consequence from (\ref{eq:4.17}) and (\ref{eq:4.14}), we get  $c\geq0 $.\\
  Then, passing to the limit in (\ref{eq:4.8}) and using (\ref{eq:4.13}),we obtain  that $u$ is a weak solution of the problem (\ref{eq:1.1}) that is $$\int_{B}\big(w(x)~|\Delta u|^{\frac{N}{2}-2}~\Delta u~~\Delta\varphi+|\nabla u|^{\frac{N}{2}-2} \nabla u.\nabla \varphi +V(x)|u|^{\frac{N}{2}-2}u\varphi\big)~dx=\int_{B}f(x,u)~ \varphi~dx,~~~~~~\mbox{for all }~~\varphi \in E.$$ Taking $\varphi=u$ as a test function, we get
$$\|u\|^{\frac{N}{2}}=\int_{B}w(x) |\Delta u|^{\frac{N}{2}} dx+\int_{B}|\nabla u|^{\frac{N}{2}}dx+\int_{B}V(x)|u|^{\frac{N}{2}}dx=
\int_{B}f(x,u)u dx\geq \frac{N}{2}\int_{B}F(x,u)dx\cdot$$
~Hence $\mathcal{J}(u)\geq 0$ . We also have by Fatou's lemma and (\ref{eq:4.14}) that $$0\leq \mathcal{J}(u)\leq \frac{2}{N}\liminf_{n\rightarrow\infty} \|u_{n}\|^{\frac{N}{2}}-\int_{B}F(x,u)dx=c.$$

So, we will finish the proof by considering  three cases for the level $c$.\\\\
{\it \underline{Case 1}.} $c=0$.~~ In this case
 $$0\leq\mathcal{ J}(u)\leq \liminf_{n\rightarrow+\infty}\mathcal{J}(u_{n})=0.$$
So, $$\mathcal{J}(u)=0$$
   and then by (\ref{eq:4.14})
$$\displaystyle\lim_{n\rightarrow +\infty}\frac{2}{N}\|u_{n}\|^{\frac{N}{2}}=\int_{B}F(x,u)dx=\frac{2}{N}\|u\|^{\frac{N}{2}}.$$
 By Brezis-Lieb's Lemma \cite{Br}, it follows that $u_{n}\rightarrow u~~\mbox{in}~~E$.\\

\noindent {\it\underline{ Case 2}.} $c>0$ and $u=0$. We prove that this case cannot happen.\\
From (\ref{eq:4.7}) and (\ref{eq:4.8}) with $v=u_{n}$, we have
$$\displaystyle\lim_{n\rightarrow +\infty}\|u_{n}\|^{\frac{N}{2}}=\frac{N}{2}c~~ \mbox{and}~~\displaystyle\lim_{n\rightarrow +\infty}\int_{B}f(x,u_{n})u_{n}dx=\frac{N}{2}c.$$
Again, by (\ref{eq:4.8}) we have
$$\big| \|u_{n}\|^{\frac{N}{2}}-\int_{B}f(x,u_{n})u_{n}dx\big|\leq C\varepsilon_{n}.$$
First we claim that there exists $q>1$ such that
\begin{equation}\label{eq:4.18}
\int_{B}|f(x,u_{n})|^{q}dx\leq C,
\end{equation}
so
$$\|u_{n}\|^{\frac{N}{2}}\leq C\varepsilon_{n}+\big(\int_{B}|f(x,u_{n})|^{q}\big)^{\frac{1}{q}}dx(\int_{B}|u_{n}|^{q'}\big)^{\frac{1}{q'}}$$
where $q'$ the conjugate of $q$. Since $(u_{n})$ converge to $u=0$ in $L^{q'}(B)$,
$$\displaystyle\lim_{n\rightarrow +\infty}\|u_{n}\|^{\frac{N}{2}}=0$$
which in contradiction with $c>0$.\\
For the proof of the claim, since $f$ has critical growth, for every
$\varepsilon>0$ and $q>1$ there exists $t_{\varepsilon}>0$ and $C>0$
such that for all $|t|\geq
t_{\varepsilon}$, we have
\begin{equation}\label{eq:4.19}
|f(x,t)|^{q}\leq
Ce^{\alpha_{0}(1+\varepsilon) t^{\gamma}}.
\end{equation}
Consequently,
$$\begin{array}{rcll}
\displaystyle\int_{B}|f(x,u_{n})|^{q}dx&=&\displaystyle\int_{\{|u_{n}|\leq
t_{\varepsilon}\}}|f(x,u_{n})|^{q}dx
+\int_{\{|u_{n}|>t_{\varepsilon}\}}|f(x,u_{n})|^{q}dx\\
 &\leq &\displaystyle\omega_{N-1} \max_{B\times [-t_{\varepsilon},t_{\varepsilon}]}|f(x,t)|^{q}+ C\int_{B}e^{\alpha_{0}(1+\varepsilon)|u_{n}|^{\gamma}}\big)dx,
 \end{array}$$
  with $\omega_{N-1}$ is the area of the unit sphere $S^{N-1}$.
Since $(\frac{N}{2}c)^{\frac{2\gamma}{N}}<\displaystyle(\frac{\alpha_{\beta}}{\alpha_{0}})$, there exists
$\eta\in(0,\frac{1}{2})$ such that
$(\frac{N}{2}c)^{\frac{2\gamma}{N}}=\displaystyle(1-2\eta)\displaystyle(\frac{\alpha_{\beta}}{\alpha_{0}})$.
On the other hand, $\|u_{n}\|^{\gamma}\rightarrow
(\frac{N}{2}c)^{\frac{2\gamma}{N}}$, so there exists $n_{\eta}>0$ such that for
all $n\geq n_{\eta}$, we get $\|u_{n}\|^{\gamma}\leq
(1-\eta)\frac{\alpha_{\beta}}{\alpha_{0}}$.
Therefore,
$$\alpha_{0}(1+\varepsilon)(\frac{|u_{n}|}{\|u_{n}\|})^{\gamma}\|u_{n}\|^{\gamma}\leq
(1+\varepsilon)(1-\eta)~\alpha_{\beta}\cdot$$
 We choose $\varepsilon >0$ small enough to get
 $$\alpha_{0}(1+\varepsilon) \|u_{n}\|^{\gamma}\leq \alpha_{\beta}\cdot $$
 Therefore, the second integral is
uniformly bounded in view of (\ref{eq:1.4}) and the claim is
proved.\\

{\it \underline{Case 3}.} $c>0$ and $u\neq 0$.~~ In this case,
we claim that $\mathcal{J}(u)=c$ and therefore, we get
$$\lim_{n\rightarrow
+\infty}\|u_{n}\|^{\frac{N}{2}}=\frac{N}{2}\big(c+\int_{B}F(x,u)dx\big)=\big(\mathcal{J}(u)+\int_{B}F(x,u)dx\big)=\|u\|^{\frac{N}{2}}.$$
Do not forgot that
 $$\mathcal{J}(u)\leq \frac{2}{N}\liminf_{n\rightarrow+\infty} \|u_{n}\|^{\frac{N}{2}}-\int_{B}F(x,u)dx=c.$$
Suppose that $\mathcal{J}(u)<c$. Then,
\begin{equation}\label{eq:4.20}
\|u\|^{\frac{N}{2}}<(\frac{N}{2}\big(c+\int_{B}F(x,u)dx\big)\big)^{\frac{2}{N}}.
\end{equation}
Set
$$v_{n}=\displaystyle\frac{u_{n}}{\|u_{n}\|}$$ and
 $$v=\displaystyle\frac{u}{(\frac{N}{2}\big(c+\displaystyle\int_{B}F(x,u)dx\big))^{\frac{2}{N}}}\cdot$$
We have $\|v_{n}\|=1$, $v_{n}\rightharpoonup v$  in $E$, $v\not\equiv 0$ and $\|v\|<1$. So, by Lemma \ref{Lionstype}, we get
$$\displaystyle \sup_{n}\int_{B}e^{p \alpha_{\beta}|v_{n}|^{\gamma}}dx<\infty,$$ for
$1<p<U(v)=(1-\|v\|^{\frac{N}{2}})^{\frac{-2\gamma}{N}}$.\\
 As in the case $(2)$, we are going to estimate $\displaystyle\int_{B}|f(x,u_{n})|^{q}dx$.\\
For  $\varepsilon>0$, one has
$$\begin{array}{rclll}
\displaystyle\int_{B}|f(x,u_{n})|^{q}dx&=&\displaystyle\int_{\{|u_{n}|\leq
t_{\varepsilon}\}}|f(x,u_{n})|^{q}dx
+\int_{\{|u_{n}|>t_{\varepsilon}\}}|f(x,u_{n})|^{q}dx\\
 &\leq &\displaystyle NV_{N} \max_{B\times [-t_{\varepsilon},t_{\varepsilon}]}|f(x,t)|^{q}+ C\int_{B}e^{\alpha_{0}(1+\varepsilon )|u_{n}|^{\gamma}}dx\\
&\leq & C_{\varepsilon}+
 C\displaystyle\int_{B}e^{\alpha_{0}(1+\varepsilon)\|u_{n}\|^{\gamma}|v_{n}|^{\gamma}}\big)dx\leq C,
 \end{array}$$
provided that $\alpha_{0}(1+\varepsilon)\|u_{n}\|^{\gamma}\leq p~~
\alpha_{\beta}$ for some $p$ such that
$1<p<U(v)=(1-\|v\|^{\frac{N}{2}})^{\frac{-2\gamma}{N}}$.\\
In fact, note that by the definition of $v$, we have
$$\displaystyle(1-\|v\|^{\frac{N}{2}})^{\frac{-2\gamma}{N}}=\displaystyle\big(\frac{(\frac{N}{2}(c+\int_{B}F(x,u)dx)}{(\frac{N}{2}(c+\int_{B}F(x,u)dx))-\|u\|^{\frac{N}{2}})}\big)^{\frac{2\gamma}{N}}=
\big(\frac{c+\int_{B}F(x,u)dx}{c-\mathcal{J}(u)}\big)^{\frac{2\gamma}{N}}$$ Since
$$\displaystyle\lim_{n\rightarrow+\infty}\|u_{n}\|^{\gamma}=(\frac{N}{2}\big(c+\int_{B}F(x,u)dx)\big)^{\frac{2\gamma}{N}},$$
then, for large enough $n$ , we have
 $$\alpha_{0}(1+\varepsilon)\|u_{n}\|^{\gamma}\leq \alpha_{0}(1+2\varepsilon)
(\frac{N}{2}\big(c+\int_{B}F(x,u)dx\big)^{\frac{2\gamma}{N}}$$
and, to get the required estimate, we just need to show that we can choose $\varepsilon> 0$
sufficiently small such that

 $$ \frac{\alpha_{0}}{\alpha_{\beta}}(1+2 \varepsilon)
<\displaystyle\big(\frac{1}{\frac{N}{2}(c-\mathcal{J}(u))}\big)^{\frac{2\gamma}{N}},$$
that is,
 $$
 (1+2\varepsilon)\displaystyle\big(\frac{N}{2}(c-\mathcal{J}(u)\big)^{\frac{2\gamma}{N}}
<\frac{\alpha_{\beta}}{\alpha_{0}}\cdot
$$
But $\mathcal{J}(u)\geq 0$ and $c<\displaystyle\frac{2}{N}(\frac{\alpha_{\beta}}{\alpha_{0}})^{\frac{N}{2\gamma}}$, then it is in fact possible to choose $\varepsilon> 0$ such that the last inequality remains valid.\\
 So, the sequence $(f(x,u_{n}))$ is bounded in $L^{q}(B)$, $q>1$.\\
  From  (\ref{eq:4.8}) with $\varphi=u_{n}-u$, we get
 \begin{align}\label{eq:4.21}
 \nonumber\int_{B}\big(w(x)~|\Delta u_{n}|^{\frac{N}{2}-2}~\Delta u_{n}~~\Delta(u_{n}-u) dx+\int_{B}|\nabla u_{n}|^{\frac{N}{2}-2}\nabla u_{n}.(\nabla u_{n}-\nabla u) dx\\+\int_{B}V(x)|u_{n}|^{\frac{N}{2}-2}u_{n}(u_{n}-u)dx
-\int_{B}f(x,u_{n})(u_{n}-u) dx=o_{n}(1).
\end{align}
On the other hand, since $u_{n}\rightharpoonup u$ weakly in $E$ then
 \begin{equation}\label{eq:4.22}
\int_{B}w(x)|\Delta u|^{\frac{N}{2}-2}\Delta u (\Delta u_{n}-\Delta u) dx+\int_{B}|\nabla u|^{\frac{N}{2}-2}\nabla u.(\nabla u_{n}-\nabla u) dx+\int_{B}V(x)(x)|u|^{\frac{N}{2}-2}u(u_{n}-u)dx
=o_{n}(1).
\end{equation}
Combining (\ref{eq:4.21}) and (\ref{eq:4.22}),  we obtain
 \begin{align}\label{eq:4.23}\nonumber &\int_{B}w(x)\big(|\Delta u_{n}|^{\frac{N}{2}-2}\Delta u_{n}-|\Delta u|^{\frac{N}{2}-2}\Delta u\big).(\Delta u_{n}-\Delta u) dx\\ \nonumber&+\int_{B}\big(|\nabla u_{n}|^{\frac{N}{2}-2}\nabla u_{n}-|\nabla u|^{\frac{N}{2}-2}\nabla u\big).(\nabla u_{n}-\nabla u) dx \\&+\int_{B}V(x)\big(|u_{n}|^{\frac{N}{2}-2}u_{n}-|u|^{\frac{N}{2}-2}u\big)(u_{n}-u)dx
-\int_{B}f(x,u_{n})(u_{n}-u) dx=o_{n}(1).\end{align}
On one hand, by the H\"{o}lder inequality, Sobolev embedding  and using the boundedness of
$\{f(x,u_{n})\}$ in $L^{q}(B)$ for $q>1$, we obtain \begin{align}\label{eq:4.24}\displaystyle\int_{B}f(x,u_{n})(u_{n}-u) dx
\leq \displaystyle\big(\int_{B}|f(x,u_{n})|^{q}\big)^{\frac{1}{q}}(\int_{B}|u_{n}-u|^{q'})^{\frac{1}{q'}}dx\rightarrow 0~~\mbox{as}~~n\rightarrow0\cdot\end{align}
 Also, using the well known inequality
\begin{equation}\label{eq:4.25}
(|x|^{p-2}x-|y|^{p-2}y).(x-y)~ \geq 2^{2-p}|x-y|^{p},~~~\forall ~~x,y \in \mathbb{R}^{m},~~ m\in\mathbb{N}~~\mbox{and}~~~~p\geq 2,
\end{equation} we get,
\begin{align} \label{4.26}\displaystyle  \int_{B}\big(|\nabla u_{n}|^{\frac{N}{2}-2}\nabla u_{n}-|\nabla u|^{\frac{N}{2}-2}\nabla u\big).(\nabla u_{n}-\nabla u) dx\geq  2^{\frac{4-N}{2}}g_{0}\displaystyle \int_{B}|\nabla u_{n}-\nabla u|^{\frac{N}{2}}dx \geq 0, \end{align}
and
\begin{align}\label{eq:4.27}&\nonumber \displaystyle\int_{B}w(x)\big(|\Delta u_{n}|^{\frac{N}{2}-2}\Delta u_{n}-|\Delta u|^{\frac{N}{2}-2}\Delta u\big).(\Delta u_{n}-\Delta u) dx+\int_{B}V(x)\big(|u_{n}|^{\frac{N}{2}-2}u_{n}-|u|^{\frac{N}{2}-2}u\big)(u_{n}-u)dx\\
 &\geq 2^{2-\frac{N}{2}}\bigg(\int_{B}w(x)\big|\Delta u_{n}-\Delta u\big|^{\frac{N}{2}} dx+\int_{B}V(x)|u_{n}-u|^{\frac{N}{2}}dx\bigg)\geq 0.\end{align}
Then,
\begin{equation}\label{eq:3.18}
0\leq  2^{2-\frac{N}{2}}\lim_{n\rightarrow+\infty}\big(\int_{B}w(x)|\Delta u_{n}-\Delta u|^{\frac{N}{2}}dx+\int_{B}|\nabla u_{n}-\nabla u|^{\frac{N}{2}}dx+\int_{B}V(x)|u_{n}-u|^{\frac{N}{2}}dx\big)=0.
\end{equation}

So,
$$ \|u_{n}-u\| \rightarrow 0~~~~ \mbox{as}~~n\rightarrow\infty.$$
By Brezis-Lieb's lemma, up to subsequence, we get
$$\displaystyle\lim_{n\rightarrow+\infty}\|u_{n}\|^{\frac{N}{2}}=\frac{N}{2}(c+\int_{B}F(x,u)dx)=\|u\|^{\frac{N}{2}}$$
and this contradicts (\ref{eq:4.20}). So, $\mathcal{J}(u)=c$ and  consequently, $u_{n} \rightarrow u~~ \mbox{in}~~E$.

\section{Proof of the main results}
In the sequel, we will estimate the minimax level of the energy $\mathcal{J}$.
We will  prove that the mountain pass level $c$ satisfies
$$c<\displaystyle\frac{2}{N}(\frac{\alpha_{\beta}}{\alpha_{0}})^{\frac{N}{2\gamma}}\cdot$$
 To this end , we will prove that there exists $z\in E$, $\|z\|=1$ such that
  \begin{equation}\label{eq:5.1}
 \max_{t\geq0}\mathcal{J}(tz)<\displaystyle\frac{2}{N}(\frac{\alpha_{\beta}}{\alpha_{0}})^{\frac{N}{2\gamma}}\cdot
  \end{equation}
  \subsection{Adams' functions}
Next, we'll utilize specific functions, known as the Adams' functions, as described in \cite{ZZ}.
  We consider the sequence defined for all $n\geq 8$  by
  \begin{equation}\label{eq:5.2}w_{n}(x)=C(N,\beta)\displaystyle \left\{
      \begin{array}{rclll}
&\displaystyle\bigg(\frac{\log (e\sqrt[N]{n} )}{\alpha_{\beta}}\bigg)^{\frac{1}{\gamma}}-\frac{|x|^{2(1-\beta)}}{2\big(\alpha_{\beta}\big)^{\frac{1}{\gamma}}
\big(\log  (e\sqrt[N]{n} )\big)^{\frac{\gamma-1}{\gamma}}}\\&+\frac{1}{2\big(\alpha_{\beta}\big)^{\frac{1}{\gamma}}(\frac{1}{n})^{\frac{2(1-\beta)}{N}}
\big(\log  (e\sqrt[N]{n} )\big)^{\frac{\gamma-1}{\gamma}}}& \mbox{ if } 0\leq |x|\leq \frac{1}{\sqrt[N]{n}}\\\\
       &\displaystyle  \frac{1}{\alpha^{\frac{1}{\gamma}}_{\beta}\big(\log e\sqrt[N]{n} \big)^{\frac{2(1-\beta)}{N}}}\bigg(\log(\frac{e}{|x|}\bigg)^{1-\beta} & \mbox{ if } \frac{1}{\sqrt[N]{n}}\leq|x|\leq \frac{1}{2}\\
       &\zeta_{n}&\mbox{ if}~~\frac{1}{2}\leq |x|\leq 1
 \end{array}
    \right.
  \end{equation}
     where $\displaystyle C(N,\beta)=\frac{(\frac{1}{2})^{\frac{2}{N}}\alpha^{\frac{1}{\gamma}}_{\beta}}{V^{\frac{2}{N}}_{N}(1-\beta)^{1-\frac{2}{N}}(N-2)} $, $\zeta_{n}\in C^{\infty}_{0,rad}(B)$ is such that\\
      $\displaystyle\zeta_{n}(\frac{1}{2})= \frac{1}{\alpha^{\frac{1}{\gamma}}_{\beta}\big(\log e\sqrt[N]{n} \big)^{\frac{2(1-\beta)}{N}}}\big(\log 2e \big)^{1-\beta}$,
      $\displaystyle\frac{\partial \zeta_{n}}{\partial r}(\frac{1}{2}) = \frac{-2(1-\beta)}{\alpha^{\frac{1}{\gamma}}_{\beta}\big(\log e\sqrt[N]{n} \big)^{\frac{2(1-\beta)}{N}}} \big(\log (2e)\big)^{-\beta}$ \newline  $\displaystyle\zeta_{n}(1)=\frac{\partial \zeta_n}{\partial r}(1)=0$ and $\xi_{n}$, $\nabla \xi_{n}$, $\Delta \xi_{n}$ are all $\displaystyle o\bigg(\frac{1}{[\log (e\sqrt[N]{n})]^{\frac{1}{\gamma}}}\bigg)$. Here,  $\displaystyle\frac{\partial \zeta_{n}}{\partial r}$ denotes the first derivative of $\zeta_{n}$ in the radial variable $r=|x|$. \\\\

 We compute $\Delta w_{n}(x)$, we get  \begin{equation*}\Delta w_{n}(x)=C(N,\beta)\displaystyle \left\{
      \begin{array}{rclll}
\displaystyle\frac{-2(1-\beta)(N-2\beta)|x|^{-2\beta}}{\alpha^{\frac{1}{\gamma}}_{\beta}\big(\log e\sqrt[N]{n} \big)^{\frac{2(1-\beta)}{N}}\big(\log  (e\sqrt[N]{n})\big)^{\frac{\gamma-1}{\gamma}}} & \mbox{ if } 0\leq |x|\leq \frac{1}{\sqrt[N]{n}}\\\\
       \displaystyle \frac{-(1-\beta)}{|x|^{2}} \frac{\bigg(\log(\frac{e}{|x|})\bigg)^{-\beta}\bigg((N-2)+\beta\big(\log \frac{e}{|x|}\big)^{-1}\bigg)}{\alpha^{\frac{1}{\gamma}}_{\beta}\big(\log e\sqrt[N]{n} \big)^{\frac{2(1-\beta)}{N}}} & \mbox{ if } \frac{1}{\sqrt[4]{n}}\leq|x|\leq \frac{1}{2}\\
       \Delta \zeta_{n}&\mbox{ if}\frac{1}{2}\leq |x|\leq 1
 \end{array}
    \right.
  \end{equation*}
  So, \begin{align}\nonumber &\frac{1}{C^{\frac{N}{2}}(N,\beta)} \|\Delta w_{n}\|_{\frac{N}{2},w}^{\frac{N}{2}}=\underbrace{NV_{N}\int^{\frac{1}{\sqrt[N]{n}}}_{0}r^{N-1}|\Delta w_{n}(x)|^{\frac{N}{2}}\big(\log \frac{e}{r}\big)^{\beta(\frac{N}{2}-1)}dr}_{I_{1}}\\ \nonumber&+\underbrace{NV_{N}\int^{\frac{1}{2}}_{\frac{1}{\sqrt[N]{n}}}r^{N-1}|\Delta w_{n}(x)|^{\frac{N}{2}}\big(\log \frac{e}{r}\big)^{\beta(\frac{N}{2}-1)}dr}_{I_{2}}\\\nonumber &+\underbrace{NV_{N}\int^{1}_{\frac{1}{2}}r^{N-1}|\Delta w_{n}(x)|^{\frac{N}{2}}\big(\log \frac{e}{r}\big)^{\beta(\frac{N}{2}-1)}dr}_{I_{3}}.\end{align}

Using integration by parts, we get,$$\begin{array}{rclll}\displaystyle I_{1}&= & \displaystyle NV_{N}\frac{2^{\frac{N}{2}}(1-\beta)^{\frac{N}{2}}(N-2\beta)^{\frac{N}{2}}}{\big(\alpha_{\beta}\big)^{\frac{N}{2\gamma}}(\log  (e\sqrt[N]{n})\big)^{1-\beta}\big(\log  (e\sqrt[N]{n})\big)^{\frac{N(\gamma-1)}{2\gamma}}}\int^{\frac{1}{\sqrt[N]{n}}}_{0} r^{N(1-\beta)-1}\big(\log \frac{e}{r}\big)^{\beta(\frac{N}{2}-1)}dr\\
&=&\displaystyle NV_{N}\frac{2^{\frac{N}{2}}(1-\beta)^{\frac{N}{2}}(N-2\beta)^{\frac{N}{2}}}{\big(\alpha_{\beta}\big)^{\frac{N}{2\gamma}}(\log  (e\sqrt[N]{n})\big)^{1-\beta}\big(\log  (e\sqrt[N]{n})\big)^{\frac{N(\gamma-1)}{2\gamma}}}\left[  \frac{r^{N(1-\beta)}}{N(1-\beta)}(\log \frac{e}{r}\big)^{\beta(\frac{N}{2}-1)} \right]^{\frac{1}{\sqrt[N]{n}}}_{0}\\&+& NV_{N}\frac{\beta(\frac{N}{2}-1)2^{\frac{N}{2}}(1-\beta)^{\frac{N}{2}}(N-2\beta)^{\frac{N}{2}}}{\big(\alpha_{\beta}\big)^{\frac{N}{2\gamma}}(\log  (e\sqrt[N]{n})\big)^{1-\beta}\big(\log  (e\sqrt[N]{n})\big)^{\frac{N(\gamma-1)}{2\gamma}}}\displaystyle\int^{\frac{1}{\sqrt[N]{n}}}_{0}\frac{r^{N(1-\beta)}-1}{N(1-\beta)} \big(\log \frac{e}{r}\big)^{\beta(\frac{N}{2}-1)-1} dr\\
&=&\displaystyle o\big(\frac{1}{\log e\sqrt[N]{n}}\big)\cdot\\

\end{array}$$
Also, \begin{align}\nonumber\displaystyle \frac{2}{N}C^{\frac{N}{2}}(N,\beta)I_{2}=\\\nonumber= &\displaystyle C^{\frac{N}{2}}(N,\beta) NV_{N}\frac{(1-\beta)^{\frac{N}{2}}}{\big(\alpha_{\beta}\big)^{\frac{N}{2\gamma}}\big(\log  (e\sqrt[N]{n})\big)^{1-\beta}}\int_{\frac{1}{\sqrt[N]{n}}}^{\frac{1}{\frac{1}{2}} }\frac{1}{r}\big(\log \frac{e}{r}\big)^{-\beta}\big ((N-2)+\beta \big( \log \frac{e}{r}\big)^{-1}\big)^{\frac{N}{2}}dr\\\nonumber
=& \nonumber\displaystyle C^{\frac{N}{2}}(N,\beta)NV_{N}\frac{(1-\beta)^{\frac{N}{2}}(N-2)^{\frac{N}{2}}}{\big(\alpha_{\beta}\big)^{\frac{N}{2\gamma}}\big(\log  (e\sqrt[N]{n})\big)^{1-\beta}}\int_{\frac{1}{\sqrt[N]{n}}}^{\frac{1}{\frac{1}{2}} }\frac{1}{r}\big(\log \frac{e}{r}\big)^{-\beta}\big (1+o( \log \frac{e}{r}\big)^{-1}\big)^{\frac{N}{2}}dr\\=& \nonumber\displaystyle C^{\frac{N}{2}}(N,\beta) NV_{N}\frac{(1-\beta)^{\frac{N}{2}}(N-2)^{\frac{N}{2}}}{\big(\alpha_{\beta}\big)^{\frac{N}{2\gamma}}\big(\log  (e\sqrt[N]{n})\big)^{1-\beta}}\bigg(\int_{\frac{1}{\sqrt[N]{n}}}^{\frac{1}{\frac{1}{2}} }\frac{1}{r}\big(\log \frac{e}{r}\big)^{-\beta}dr +\int_{\frac{1}{\sqrt[N]{n}}}^{\frac{1}{\frac{1}{2}} }\frac{1}{r}\big(\log \frac{e}{r}\big)^{-\beta}o( \log \frac{e}{r}\big)^{-1}dr\bigg)\\=& \nonumber\displaystyle C^{\frac{N}{2}}(N,\beta)NV_{N}\frac{(1-\beta)^{\frac{N}{2}}(N-2)^{\frac{N}{2}}}{\big(\alpha_{\beta}\big)^{\frac{N}{2\gamma}}\big(\log  (e\sqrt[N]{n})\big)^{1-\beta}}\left [ \frac{1}{1-\beta}\big( \log \frac{e}{r}\big)^{1-\beta}\right]^{\frac{1}{\sqrt[N]{n}}}_{\frac{1}{2}}\\ \nonumber&-
C^{\frac{N}{2}}(N,\beta)NV_{N}\frac{(1-\beta)^{\frac{N}{2}}(N-2)^{\frac{N}{2}}}{\big(\alpha_{\beta}\big)^{\frac{N}{2\gamma}}\big(\log  (e\sqrt[N]{n})\big)^{1-\beta}}\bigg(\int_{\frac{1}{\sqrt[N]{n}}}^{\frac{1}{\frac{1}{2}} } \frac{1}{r}\big(\log \frac{e}{r}\big)^{-\beta}o( \log \frac{e}{r}\big)^{-1}dr\bigg)\\\nonumber
=& \nonumber\displaystyle 1+ o\big(\frac{1}{(\log e\sqrt[N]{n})^{1-\beta}}\big)\cdot
\end{align}
and $I_{3}=\displaystyle  o\big(\frac{1}{(\log e\sqrt[N]{n})^{\frac{2}{\gamma}}}\big).$ Then $\frac{2}{N}\|\Delta w_{n}\|_{\frac{N}{2},w}^{\frac{N}{2}}=1+o\big(\frac{1}{(\log e\sqrt[N]{n})^{\frac{N}{\gamma}}}\big)$.
\subsection{Key lemmas}

   \begin{lemma}\label{lem6} Assume V(x) is continuous and $(V_{1})$ is satisfied. Then there holds
    $\displaystyle \lim_{n\rightarrow+\infty}\|w_{n}\|^{\frac{N}{2}}=1.$

\end{lemma}
    {\it Proof}~~We have
    \begin{align}\nonumber\|w_{n}\|^{\frac{N}{2}}&=  \displaystyle\frac{2}{N}\int_{B} w(x)|\Delta w_{n}|^{\frac{N}{2}}dx+\frac{2}{N}\int_{B}|\nabla w_{n}|^{\frac{N}{2}}dx+\frac{2}{N}\int_{B}V~~|w_{n}|^{\frac{N}{2}}dx \\\nonumber
&=\displaystyle 1+o\big(\frac{1}{(\log e\sqrt[N]{n})^{\frac{2}{\gamma}}}\big)+\nonumber\frac{2}{N}\int_{0\leq |x|\leq \frac{1}{\sqrt[N]{n}}}V(x)|w|^{\frac{N}{2}}_{n}dx+\frac{2}{N}\int_{\sqrt[N]{n}\leq |x|\leq \frac{1}{2}}V(x) |w_{n}|^{\frac{N}{2}}dx\\ \nonumber &+\displaystyle\frac{2}{N}\int_{|x|\geq \frac{1}{2}}V(x)|\zeta|^{\frac{N}{2}}_{n}dx\\&\nonumber+\displaystyle\underbrace{\frac{2}{N}\int_{0\leq |x|\leq \frac{1}{\sqrt[N]{n}}}| \nabla w_{n}|^{\frac{N}{2}}dx}_{I_{1}'}+\underbrace{\frac{2}{N}\int_{\sqrt[N]{n}\leq |x|\leq \frac{1}{2}}| \nabla w_{n}|^{\frac{N}{2}}dx}_{I_{2}'}+\underbrace{\frac{2}{N}\int_{|x|\geq \frac{1}{2}}|\nabla \zeta_{n}|^{\frac{N}{2}}dx}_{I_{3}'}\cdot
\end{align}
We have,$$\begin{array}{rclll}\displaystyle I'_{1}&= & \displaystyle 2V_{N}\frac{C^{\frac{N}{2}}(N,\beta)(1-\beta)^{\frac{N}{2}}}{\alpha^{\frac{N}{2\gamma}}_{\beta}\big(\log  (\sqrt[N]{n})\big)^{\frac{N(\gamma-1)}{2\gamma}}}\int^{\frac{1}{\sqrt[N]{n}}}_{0} r^{N(2-\beta)-1}dr\\
&=&\displaystyle 2V_{N}\frac{C^{\frac{N}{2}}(N,\beta)(1-\beta)^{\frac{N}{2}}}{\alpha^{\frac{N}{2\gamma}}_{\beta}\big(\log  (\sqrt[N]{n})\big)^{\frac{N(\gamma-1)}{2\gamma}}}\left[  \frac{r^{N(2-\beta)}}{N(2-\beta)} \right]^{\frac{1}{\sqrt[N]{n}}}_{0}
\\&=&\displaystyle 2V_{N}\frac{C^{\frac{N}{2}}(N,\beta)(1-\beta)^{\frac{N}{2}}}{\alpha^{\frac{N}{2\gamma}}_{\beta}n^{2-\beta}N(2-\beta)\log  (\sqrt[N]{n})\big)^{\frac{N(\gamma-1)}{2\gamma}}}\\
&=&\displaystyle o\big(\frac{1}{n^{2-\beta}\log e\sqrt[N]{n}}\big)\cdot

\end{array}$$
Also, using the fact that the function $r\mapsto r^{\frac{N}{2}-1}\big(\log \frac{e}{r}\big)^{-\frac{N}{2}\beta}$ is increasing on $[0,1]$, we get  $$\begin{array}{rclll}\displaystyle I'_{
2}&= & \displaystyle 2V_{N}C^{\frac{N}{2}}(N,\beta)\frac{(1-\beta)^{\frac{N}{2}}}{\big(\alpha_{\beta}\big)^{\frac{N}{2\gamma}}\big(\log  (e\sqrt[N]{n})\big)^{1-\beta}}\int_{\frac{1}{\sqrt[N]{n}}}^{\frac{1}{2}} r^{\frac{N}{2}-1}\big(\log \frac{e}{r}\big)^{-\frac{N}{2}\beta}dr\\
&\leq & \displaystyle 2V_{N}C^{\frac{N}{2}}(N,\beta)\frac{(1-\beta)^{\frac{N}{2}}}{\big(\alpha_{\beta}\big)^{\frac{N}{2\gamma}}\big(\log  (e\sqrt[N]{n})\big)^{1-\beta}}(\frac{1}{2})^{\frac{N}{2}-1}\big(\log 2e\big)^{-\frac{N}{2}\beta}\\
&=&\displaystyle o\bigg(\frac{1}{[\log (e\sqrt[N]{n})]^{1-\beta}}\bigg)\\
\end{array}$$
and $I'_{3}=\displaystyle  o\big(\frac{1}{(\log e\sqrt[N]{n})^{\frac{2}{\gamma}}}\big).$
    For $\displaystyle|x|\leq \frac{1}{\sqrt[N]{n}}$, $$\displaystyle |w_{n}|^{\frac{N}{2}}\leq C(N,\beta)^{\frac{N}{2}} \bigg(\frac{\log (e\sqrt[N]{n} )}{\alpha_{\beta}}\bigg)^{\frac{1}{\gamma}}+ \frac{1}{2\big(\alpha_{\beta}\big)^{\frac{1}{\gamma}}(\frac{1}{n})^{\frac{2(1-\beta)}{N}}
\big(\log  (e\sqrt[N]{n} )\big)^{\frac{\gamma-1}{\gamma}}}\bigg)^{\frac{N}{2}}\cdot$$Then,
  \begin{align}\nonumber&\int_{0\leq |x|\leq \frac{1}{\sqrt[N]{n}}}V(x)|w_{n}|^{\frac{N}{2}}dx\leq\\\nonumber & NV_{N} mC(N,\beta)^{\frac{N}{2}}\bigg(\bigg(\frac{\log (e\sqrt[N]{n} )}{\alpha_{\beta}}\bigg)^{\frac{1}{\gamma}}+ \frac{1}{2\big(\alpha_{\beta}\big)^{\frac{1}{\gamma}}(\frac{1}{n})^{\frac{2(1-\beta)}{N}}
\big(\log  (e\sqrt[N]{n} )\big)^{\frac{\gamma-1}{\gamma}}}\bigg)^{\frac{N}{2}}\int^{\frac{1}{\sqrt[N]{n}}}_{0} r^{N-1}dr=o_{n}(1)\end{align}
 Also,$$\int_{\frac{1}{\sqrt[N]{n}}\leq |x|\leq \frac{1}{2}}V(x)|w_{n}|^{\frac{N}{2}}dx\leq \frac{NV_{N} mC(N,\beta)^{\frac{N}{2}} }{\big(\alpha_{\beta}\big)^{\frac{N}{2\gamma}}\big(\log  (e\sqrt[N]{n})\big)^{1-\beta}}\int^{\frac{1}{2}}_{\frac{1}{\sqrt[N]{n}}}r^{N-1}\big(\log(\frac{e}{r}\big)\big)^{\frac{2(1-\beta)}{N}}dr$$
 Using the fact that the function $r\mapsto r^{N-1}\big(\log \frac{e}{r}\big)^{\frac{2(1-\beta)}{N}}$ is increasing on $[0,1]$,   we obtain $$\int_{\frac{1}{\sqrt[N]{n}}\leq |x|\leq \frac{1}{2}}V(x)|w_{n}|^{\frac{N}{2}}_{n}dx\leq \frac{NV_{N} mC(N,\beta)^{\frac{N}{2}} }{\big(\alpha_{\beta}\big)^{\frac{N}{2\gamma}}\big(\log  (e\sqrt[N]{n})\big)^{1-\beta}}\frac{1}{2^{N-1}}\big(\log  (2e)\big)^{\frac{2(1-\beta)}{N}}=o_{n}(1).$$
 Finaly,$$\int_{{\frac{1}{2}\leq |x|\leq 1}} V(x) |w_{n}|^{\frac{N}{2}}dx\leq m\int_{|x|\geq \frac{1}{2}}|\zeta_{n}|^{\frac{N}{2}}dx= o_{n}(1)$$
 Then, $\|\ w_{n}\|^{\frac{N}{2}}\leq1+o\big(\frac{1}{(\log e\sqrt[N]{n})^{\frac{2}{\gamma}}}\big)$.
 \\
Now,  from the definition of $w_{n}$, it is easy to see that  $$-\frac{|x|^{2(1-\beta)}}{2\big(\alpha_{\beta}\big)^{\frac{1}{\gamma}}
\big(\log  (e\sqrt[N]{n} )\big)^{\frac{\gamma-1}{\gamma}}}+\frac{1}{2\big(\alpha_{\beta}\big)^{\frac{1}{\gamma}}(\frac{1}{n})^{\frac{2(1-\beta)}{N}}
\big(\log  (e\sqrt[N]{n} )\big)^{\frac{\gamma-1}{\gamma}}}\geq0 ~~\mbox{for all}~~ 0\leq |x|\leq \frac{1}{\sqrt[N]{n}}\cdot$$ Then, for all $0\leq |x|\leq \frac{1}{\sqrt[N]{n}}$, $\displaystyle| w_{n}|^{\frac{N}{2}}\geq \displaystyle\displaystyle\bigg(\frac{\log (e\sqrt[N]{n} )}{\alpha_{\beta}}\bigg)^{\frac{N}{2\gamma}}\cdot$  So, using the condition $(V_{1})$, we get
  $$\int_{0\leq |x|\leq\frac{1}{\sqrt[N]{n}}}V(x)|w_{n}|^{\frac{N}{2}}dx\geq V_{0} NV_{N}\displaystyle\bigg(\frac{\log (e\sqrt[N]{n} )}{\alpha_{\beta}}\bigg)^{\frac{N}{2\gamma}}\int^{\frac{1}{\sqrt[N]{n}}}_{0}r^{N-1}dr=o_{n}(1)$$
  and using the fact that the function $r\mapsto r^{N-1}\big(\log(\frac{e}{r}\big)\big)^{\frac{2(1-\beta)}{N}}$ is increasing on $[0,\frac{1}{2}]$, we get $$\int_{\frac{1}{\sqrt[4]{n}}\leq |x|\leq \frac{1}{2}}V(x)|w_{n}|^{\frac{N}{2}}dx\geq \frac{NV_{N} mC(N,\beta)^{\frac{N}{2}} }{\big(\alpha_{\beta}\big)^{\frac{N}{2\gamma}}\big(\log  (e\sqrt[N]{n})\big)^{1-\beta}}\int^{\frac{1}{2}}_{\frac{1}{\sqrt[N]{n}}}r^{N-1}\big(\log(\frac{e}{r}\big)\big)^{\frac{2(1-\beta)}{N}}dr=o_{n}(1)$$
Consequently, $ 1+o_{1}\big(\frac{1}{(\log e\sqrt[N]{n})^{\frac{N}{2\gamma}}}\big)\leq\|w_{n}\|^{\frac{N}{2}}\leq 1+o_{2}\big(\frac{1}{(\log e\sqrt[N]{n})^{\frac{N}{2\gamma}}}\big)$. The Lemma is proved. \\

  \subsection{Min-Max level estimate}
We are aiming to obtain the desired estimate.

\begin{lemma}\label{lem4.2}~~For the sequence $(v_{n})$ identified by (\ref{eq:5.2}), there exists $n\geq1$ such that
  \begin{equation}\label{eq:5.3}
 \max_{t\geq 0}\mathcal{J}(tv_{n})<\displaystyle\frac{2}{N}(\frac{\alpha_{\beta}}{\alpha_{0}})^{\frac{N}{2\gamma}}\cdot
       \end{equation}
\end{lemma}
{\it Proof}~~ By contradiction, suppose that for all $n\geq1$,
$$ \max_{t\geq 0}\mathcal{J}(tv_{n})\geq\displaystyle\frac{2}{N}(\frac{\alpha_{\beta}}{\alpha_{0}})^{\frac{N}{2\gamma}}\cdot$$
Therefore,
for any $n\geq1$, there exists $t_{n}>0$ such that
$$\max_{t\geq0}\mathcal{J}(tv_{n})=\mathcal{J}(t_{n}v_{n})\geq \displaystyle\frac{2}{N}(\frac{\alpha_{\beta}}{\alpha_{0}})^{\frac{N}{2\gamma}}$$
and so,
$$\frac{2}{N}t_{n}^{\frac{N}{2}}-\int_{B}F(x,t_{n}v_{n})dx\geq \displaystyle\frac{2}{N}(\frac{\alpha_{\beta}}{\alpha_{0}})^{\frac{N}{2\gamma}}\cdot$$
Then, by using $(H_{1})$
   \begin{equation}\label{eq:5.4}
t_{n}^{\frac{N}{2}}\geq\displaystyle(\frac{\alpha_{\beta}}{\alpha_{0}})^{\frac{N}{2\gamma}}\cdot
       \end{equation}
On the other hand,
 $$\frac{d}{dt}\mathcal{J}(tv_{n})\big|_{t=t_{n}}=t^{\frac{N}{2}-1}_{n}-\int_{B}f(x,t_{n}v_{n})v_{n}dx=0,$$
then,
\begin{equation}\label{eq:5.5}
t_{n}^{\frac{N}{2}}=\int_{B}f(x,t_{n}v_{n})t_{n}v_{n}dx.
\end{equation}

Now, we claim that the sequence $(t_{n})$ is bounded in $(0,+\infty)$. Indeed,
 it follows from $(H_{5})$ that for all $\varepsilon>0$, there exists
$t_{\varepsilon}>0$ such that
\begin{equation}\label{eq:5.6}
f(x,t)t\geq
(\gamma_{0}-\varepsilon)e^{\alpha_{0}|t|^{\gamma}}~~\forall
|t|\geq t_{\varepsilon},~~\mbox{uniformly in}~~x\in B.
\end{equation}
Using  (\ref{eq:5.5}), we get
\begin{equation*}\label{eq:5.13}
t_{n}^{\frac{N}{2}}=\int_{B}f(x,t_{n}v_{n})t_{n}v_{n}dx\geq
\int_{0\leq |x|\leq \frac{1}{\sqrt[N]{n}}}f(x,t_{n}v_{n})t_{n}v_{n}dx
\cdot\end{equation*}
We have  for all $0\leq |x|\leq \frac{1}{\sqrt[N]{n}}$, $\displaystyle w^{\gamma}_{n}\geq \displaystyle\bigg(\frac{\log (e\sqrt[N]{n} )}{\alpha_{\beta}}\bigg)\cdot$
From (\ref{eq:5.4}) and the result of Lemma \ref{lem6},  $$t_{n} v_{n}\geq\frac{t_{n}}{\|w_{n}\|}\big(\frac{\log (e\sqrt[N]{n} )}{\alpha_{\beta}}\big)^{\frac{1}{\gamma}} \rightarrow\infty~~\mbox{as}~~n\rightarrow+\infty.$$ Hence,
it follows from  (\ref{eq:5.6}) that for all $\varepsilon >0$, there exists $n_{0}$ such that for all $n\geq n_{0}$
\begin{equation*}
t^{\frac{N}{2}}_{n}\geq \displaystyle (\gamma_{0}-\varepsilon)\int_{0\leq |x|\leq \frac{1}{\sqrt[N]{n}}}e^{\alpha_{0}t^{\gamma}_{n}|v|^{\gamma}_{n}}dx
  \end{equation*}
and  \begin{equation}\label{eq:5.7}
  t^{\frac{N}{2}}_{n}\geq \displaystyle NV_{N} (\gamma_{0}-\varepsilon)\int_{0}^{\frac{1}{\sqrt[N]{n}}}r^{N-1}e^{\alpha_{0}t^{\gamma}_{n} \displaystyle\big(\frac{\log (e\sqrt[N]{n})}{\|w_{n}\|^{\gamma}\alpha_{\beta}}\big)}dr\cdot
\end{equation}
Hence,
 \begin{equation*}\label{eq:5.13}
1 \geq  \displaystyle NV_{N}  (\gamma_{0}-\varepsilon)~~\displaystyle e^{\alpha_{0}t^{\gamma}_{n} \displaystyle\big(\frac{\log (e\sqrt[N]{n})}{\|w_{n}\|^{\gamma}\alpha_{\beta}}\big)-\log Nn-\frac{N}{2}\log t_{n}}.
 \end{equation*}
 Therefore $(t_{n})$ is bounded. Also, we have from the formula  (\ref{eq:5.5}) that
$$\displaystyle\lim_{n\rightarrow+\infty} t_{n}^{\frac{N}{2}}\geq\displaystyle(\frac{\alpha_{\beta}}{\alpha_{0}})^{\frac{N}{2\gamma}}\cdot$$
Now, suppose that
$$\displaystyle\lim_{n\rightarrow+\infty} t_{n}^{\frac{N}{2}}>\displaystyle(\frac{\alpha_{\beta}}{\alpha_{0}})^{\frac{N}{2\gamma}},$$
then for $n$ large enough, there exists some $\delta>0$ such that $ t_{n}^{\gamma}\geq \frac{\alpha_{\beta}}{\alpha_{0}}+\delta$. Consequently the right hand side of (\ref{eq:5.7}) tends to infinity and this contradicts the boudness of  $(t_{n})$. Since $(t_{n})$ is bounded, we get
\begin{equation}\label{eq:5.8}
 \displaystyle
 \lim_{n\rightarrow+\infty}t_{n}^{\frac{N}{2}}=\displaystyle (\frac{\alpha_{\beta}}{\alpha_{0}})^{\frac{N}{2\gamma}}\cdot
 \end{equation}

Let consider the sets $$\mathcal{A}_{n}=\{x\in B| t_{n}v_{n}\geq
t_{\varepsilon}\}~~\mbox{and}~~\mathcal{C}_{n}=B\setminus \mathcal{A}_{n}.$$ we have,
$$\begin{array}{rclll} t_{n}^{\frac{N}{2}}&=&\displaystyle\int_{B}f(x,t_{n}v_{n})t_{n}v_{n}dx=\int_{\mathcal{A}_{n}}f(x,t_{n}v_{n})t_{n}v_{n}dx+\int_{\mathcal{C}_{n}}f(x,t_{n}v_{n})t_{n}v_{n} $$\\
&\geq& \displaystyle(\gamma_{0}-\varepsilon)\int_{\mathcal{A}_{n}}e^{\alpha_{0}t_{n}^{\gamma}v_{n}^{\gamma}}dx + \int_{\mathcal{C}_{n}}f(x,t_{n}v_{n})t_{n}v_{n}dx\\
&=&\displaystyle(\gamma_{0}-\varepsilon)\int_{B}e^{\alpha_{0}t_{n}^{\gamma}v_{n}^{\gamma}}dx-
(\gamma_{0}-\varepsilon)\int_{\mathcal{C}_{n}}e^{\alpha_{0}t_{n}^{\gamma}v_{n}^{\gamma}}dx\\ &+&\displaystyle\int_{\mathcal{C}_{n}}f(x,t_{n}v_{n})t_{n}v_{n}dx.
\end{array}$$
Since $v_{n}\rightarrow 0 ~~\mbox{a.e in }~~B$,
$\chi_{\mathcal{C}_{n}}\rightarrow1~~\mbox{a.e in}~~B$, therefore using the dominated convergence
theorem, we get $$\displaystyle\int_{\mathcal{C}_{n}}f(x,t_{n}v_{n})t_{n}v_{n}dx\rightarrow 0~~\mbox{and}~~\int_{\mathcal{C}_{n}}e^{\alpha_{0}t_{n}^{\gamma}v_{n}^{\gamma}}dx\rightarrow NV_{N}\cdot$$
Then,$$\ \lim_{n\rightarrow+\infty} t_{n}^{\frac{N}{2}}=\displaystyle(\frac{\alpha_{\beta}}{\alpha_{0}})^{\frac{N}{2\gamma}}\geq(\gamma_{0}-
\varepsilon)\lim_{n\rightarrow+\infty}\int_{B}e^{\alpha_{0}t_{n}^{\gamma}v_{n}^{\gamma}}dx-(\gamma_{0}-
\varepsilon)NV_{N}\cdot$$
On the other hand,
$$\int_{B}e^{\alpha_{0}t_{n}^{\gamma}v_{n}^{\gamma}}dx \geq \int_{\frac{1}{\sqrt[N]{n}}\leq |x|\leq \frac{1}{2}}e^{\alpha_{0}t_{n}^{\gamma}v_{n}^{\gamma}}dx+\int_{\mathcal{C}_{n}} e^{\alpha_{0}t_{n}^{\gamma}v_{n}^{\gamma}}dx\cdot$$

Then, using (\ref{eq:5.4})  $$\lim_{n\rightarrow+\infty} t_{n}^{\frac{N}{2}}\geq \lim_{n\rightarrow+\infty}\displaystyle(\gamma_{0}-\varepsilon)\int_{B}e^{\alpha_{0}t_{n}^{\gamma}v_{n}^{\gamma}}dx\geq \displaystyle\lim_{n\rightarrow+\infty}(\gamma_{0}-\varepsilon) NV_{N}\int^{\frac{1}{2}}_{\frac{1}{\sqrt[N]{n}}}r^{N-1} e^{C^{\gamma}(N,\beta)\frac{ \big(\log \frac{e}{r}\big)^{\frac{N}{N-2}}}{(\log (e\sqrt[N]{n}))^{\frac{2}{N-2}}\|w_{n}\|^{\gamma}}} dr.$$

Therefore,  making the change of variable $$s=\frac{C(N,\beta)^{\gamma}(\log \frac{e}{r})}{(\log (e\sqrt[N]{n}))^{\frac{2}{N-2}}\|w_{n}\|^{\gamma}}=P\frac{(\log \frac{e}{r})}{\|w_{n}\|^{\gamma}},~~\mbox{with}~~P=\frac{C(N,\beta)^{\gamma}}{(\log (e\sqrt[N]{n}))^{\frac{2}{N-2}}}$$
we get
$$\begin{array}{rclll}\displaystyle \lim_{n\rightarrow+\infty} t_{n}^{\frac{N}{2}}&\geq & \displaystyle\lim_{n\rightarrow+\infty}\displaystyle(\gamma_{0}-\varepsilon)\int_{B}e^{\alpha_{0}t_{n}^{\gamma}v_{n}^{\gamma}}dx\\\\
&\geq&\displaystyle \lim_{n\rightarrow+\infty} NV_{N}(\gamma_{0}-\varepsilon)\frac{\|w_{n}\|^{\gamma}}{P}\int^{\frac{P \log(e\sqrt[N]{n})}{\|w_{n}\|^{\gamma}}}_{\frac{P \log(2e)}{\|w_{n}\|^{\gamma}}}\displaystyle e^{N(1-\frac{s\|w_{n}\|^{\gamma}}{P})+\frac{\|w_{n}\|^{\frac{2\gamma}{N-2}}}{P^{\frac{N}{N-2}}}s^{\frac{N}{N-2}}}ds\\\\
&\geq&\displaystyle \lim_{n\rightarrow+\infty} NV_{N}(\gamma_{0}-\varepsilon)\frac{\|w_{n}\|^{\gamma}}{P}e^{N}\int^{\frac{P \log(e\sqrt[N]{n})}{\|w_{n}\|^{\gamma}}}_{\frac{P \log(2e)}{\|w_{n}\|^{\gamma}}}\displaystyle e^{-\frac{N}{P}\|w_{n}\|^{\gamma}s}ds\\
&=&\displaystyle\lim_{n\rightarrow+\infty}(\gamma_{0}-\varepsilon)NV_{N}\frac{e^{N}}{N}\big(-e^{-N\log (e\sqrt[N]{n})}+e^{-N\log (2e)}\big)\\
&=&\displaystyle(\gamma_{0}-\varepsilon)V_{N}e^{N(1-\log (2e))}\cdot
\end{array}$$

It follows that
\begin{equation*}\label{eq:4.15}
\displaystyle\displaystyle(\frac{\alpha_{\beta}}{\alpha_{0}})^{\frac{N}{2\gamma}}\geq (\gamma_{0}-\varepsilon)V_{N}e^{N(1-\log (2e))}
\end{equation*}
for all $\varepsilon>0$. So,
$$\gamma_{0}\leq\displaystyle \frac{(\frac{\alpha_{\beta}}{\alpha_{0}})^{\frac{N}{2\gamma}}}{V_{N}e^{N(1-\log (2e))}}, $$
   which is in contradiction with  the condition $(H_{5})$. \\Now by Proposition \ref{propPS}, the functional $\mathcal{J}$ satisfies the $(PS)_{c}$ condition at a level $c<\displaystyle\displaystyle\frac{2}{N}(\frac{\alpha_{\beta}}{\alpha_{0}})^{\frac{N}{2\gamma}}$ . Also, by Proposition \ref{prg} , we deduce that the functional $\mathcal{J}$  has a nonzero critical point $u$ in $E$. The Theorem \ref{th2} is proved.
    \section*{ Statements and Declarations:}

 We declare that this manuscript is original, has not been published before and is not currently being considered for
publication elsewhere.

We confirm that the manuscript has been read and approved and that there are no other persons who satisfied the criteria for
authorship but are not listed.
 \section*{ Competing Interests:}

The authors declare that they have no known competing financial interests or personal relationships that could have appeared to influence the work reported in this paper.

\end{document}